\documentclass[12pt,a4paper]{amsart}
\usepackage{amsfonts,amscd,amsmath,amssymb,amsthm,mathrsfs}
\usepackage{array,latexsym,float,enumerate,ifpdf}
\usepackage[usenames,dvipsnames]{color}
\usepackage{graphicx}
\usepackage[top=2.75cm, bottom=2.75cm, left=2.5cm, right=2.5cm,twoside=false]{geometry}
\usepackage[all]{xy}
\numberwithin{equation}{section}

\makeatletter
\def\subsection{\@startsection{subsection}{3}%
  \z@{.5\linespacing\@plus.7\linespacing}{.5\linespacing}%
  {\normalfont\itshape}}
\makeatother

\makeatletter
\renewenvironment{proof}[1][\proofname]{%
  \par\pushQED{\qed}\normalfont%
  \topsep6\p@\@plus6\p@\relax
  \trivlist\item[\hskip\labelsep\bfseries#1\@addpunct{.}]%
  \ignorespaces}{%
  \popQED\endtrivlist\@endpefalse}
\makeatother

\theoremstyle{plain}
\newtheorem{theorem}{Theorem}[section]
\newtheorem{lemma}[theorem]{Lemma}
\newtheorem{proposition}[theorem]{Proposition}
\newtheorem{corrolary}[theorem]{Corollary}

\newtheorem*{claim}{Claim} 
\theoremstyle{definition}
\newtheorem{definition}[theorem]{Definition}
\newtheorem{example}[theorem]{Example}
\newtheorem{noname}[theorem]{}
\newtheorem{remark}[theorem]{Remark}
\newtheorem{construction}[theorem]{Construction}
\newtheorem{notation}[theorem]{Notation}
\theoremstyle{remark}
\newtheorem*{smallremark}{Remark}

\newtheorem{case}{Case} \makeatletter \@addtoreset{case}{theorem}\makeatother

\newcommand{\bthm}{\begin{theorem}}
\newcommand{\bprop}{\begin{proposition}}
\newcommand{\blem}{\begin{lemma}}
\newcommand{\bcor}{\begin{corrolary}}
\newcommand{\brem}{\begin{remark}}
\newcommand{\bdfn}{\begin{definition}}
\newcommand{\bitem}{\begin{itemize}}
\newcommand{\bex}{\begin{example}}
\newcommand{\bno}{\begin{noname}}
\newcommand{\bsrem}{\begin{smallremark}}
\newcommand{\bnot}{\begin{notation}}
\newcommand{\bcon}{\begin{construction}}
\newcommand{\bca}{\begin{case}}
\newcommand{\bcl}{\begin{claim}}
\newcommand{\beq}{\begin{equation}}

\newcommand{\eeq}{\end{equation}}
\newcommand{\ecl}{\end{claim}}
\newcommand{\eca}{\end{case}}
\newcommand{\econ}{\end{construction}}
\newcommand{\enot}{\end{notation}}
\newcommand{\esrem}{\end{smallremark}}
\newcommand{\eno}{\end{noname}}
\newcommand{\eex}{\end{example}}
\newcommand{\eitem}{\end{itemize}}
\newcommand{\ethm}{\end{theorem}}
\newcommand{\eprop}{\end{proposition}}
\newcommand{\elem}{\end{lemma}}
\newcommand{\ecor}{\end{corrolary}}
\newcommand{\erem}{\end{remark}}
\newcommand{\edfn}{\end{definition}}
\newcommand{\benum}{\begin{enumerate}}
\newcommand{\eenum}{\end{enumerate}}

\newcommand{\wt}{\widetilde}
\newcommand{\cal}[1]{\mathcal{#1}}

\newcommand{\ds}{\displaystyle}

\def\8{\infty}
\def\.{\cdot}
\def\PP{\mathbb{P}}

\def\C{\mathbb{C}}

\def\Q{\mathbb{Q}}

\def\:{\colon}
\def\map{\dashrightarrow}

\renewcommand{\iff}{\Leftrightarrow}

\def\med{\medskip}
\def\ssk{\smallskip}

\def\Bk{\operatorname{Bk}}

\def\Sing{\operatorname{Sing}}

\def\Supp{\operatorname{Supp}}
\def\Pic{\operatorname{Pic}}

\def\Exc{\operatorname{Exc}}

\def\Spec{\operatorname{Spec}}
\def\NS{\operatorname{NS}}

\def\ind{\operatorname{ind}}

\def\tip{\operatorname{tip}}

\ifpdf \usepackage[linkbordercolor={0 0 1}]{hyperref} \else \usepackage[hypertex,linkbordercolor={0 0 1}]{hyperref} \fi

\begin{document}

\title[The Coolidge-Nagata conjecture, part I]{The Coolidge-Nagata conjecture, part I}

\author[Karol Palka]{Karol Palka}
\address{Karol Palka: Institute of Mathematics, Polish Academy of Sciences, ul. \'{S}niadeckich 8, 00-656 Warsaw, Poland}
\thanks{The author was partially supported by the Foundation For Polish Science under the Homing Plus programme cofinanced from the European Union, Regional Development Fund and by the National Science Center (NCN), Poland, grant No. 2012/05/D/ST1/03227.}

\email{palka@impan.pl}

\begin{abstract}Let $E\subseteq \PP^2$ be a complex rational cuspidal curve contained in the projective plane and let $(X,D)\to (\PP^2,E)$ be the minimal log resolution of singularities. Applying the log Minimal Model Program to $(X,\frac{1}{2}D)$ we prove that if $E$ has more than two singular points or if $D$, which is a tree of rational curves, has more than six maximal twigs or if $\PP^2\setminus E$ is not of log general type then $E$ is Cremona equivalent to a line, i.e.\ the Coolidge-Nagata conjecture for $E$ holds. We show also that if $E$ is not Cremona equivalent to a line then the morphism onto the minimal model contracts at most one irreducible curve not contained in $D$.
\end{abstract}

\maketitle

\section{Main results and strategy}\label{sec:main result}

Let $\bar E\subseteq \PP^2$ be a complex planar rational curve which is \emph{cuspidal}, i.e.\ which has only locally analytically irreducible (unibranched) singularities. Equivalently, it can be defined as an image of a singular embedding of a complex projective line into a complex projective plane, i.e.\ of a morphism $\PP^1\to \PP^2$ which is 1-1 on closed points. We say that two planar curves are \emph{Cremona equivalent} if one of them is a proper transform of the other under some Cremona transformation of $\PP^2$. Not all rational curves on $\PP^2$ are Cremona equivalent to a line (a general rational curve of degree at least six is not, see \ref{rem:deg<=5_rectifiable}) and, clearly, the proper transform of $\bar E$ under a Cremona transformation does not have to be cuspidal. Therefore, the conjecture that nevertheless all cuspidal curves are Cremona equivalent to a line, which is known as the Coolidge-Nagata conjecture, comes as a surprise. It has been studied for a long time.\footnote{Coolidge \cite[Book IV,\S II.2]{Coolidge} and Nagata \cite{Nagata} studied planar rational curves and their behaviour under the action of the Cremona group. The problem of determining which rational curves are Cremona equivalent to a line is known as the  'Coolidge-Nagata problem'.} Let $c$ be the number of cusps of $\bar E$ and let $(X,D)\to (\PP^2,\bar E)$ be the minimal log resolution of singularities. In \cite{Palka-minimal_models} we proved that $$c\leq 9-2p_2(\PP^2,\bar E),$$ where $p_2(\PP^2,\bar E)=h^0(2K_X+D)$. Let $E$ be the proper transform of $\bar E$ on $X$. The Coolidge-Nagata conjecture for $\bar E\subseteq \PP^2$ is known to be equivalent to the vanishing of $h^0(2K_X+E)$, so if it fails for $\bar E$ then we get a lower bound $p_2(\PP^2,\bar E)\geq h^0(2K_X+E)\geq 1$. The higher lower bound on $p_2(\PP^2,\bar E)$ we can prove, the bigger is the restriction on $c$ (in fact also on many other parameters describing the geometry of $\bar E\subseteq \PP^2$), and hence the closer we are to proving the conjecture. Deepening the analysis of minimal models of $(X,\frac{1}{2}D)$ started in \cite{Palka-minimal_models} (which is an analog of the 'theory of peeling' \cite[\S 2.3]{Miyan-OpenSurf} for half-integral divisors) we show here the following result.

\bthm\label{thm:MAIN1} Let $\bar E\subseteq \PP^2$ be a complex rational cuspidal curve which is not Cremona equivalent to a line and let $(X,D)\to (\PP^2,\bar E)$ be the minimal log resolution of singularities. Then $p_2(\PP^2,\bar E)\in\{3,4\}$. Equivalently, $(K_X+D)^2\in\{1,2\}$. \ethm

\begin{figure}[h]\centering\includegraphics[scale=0.5]{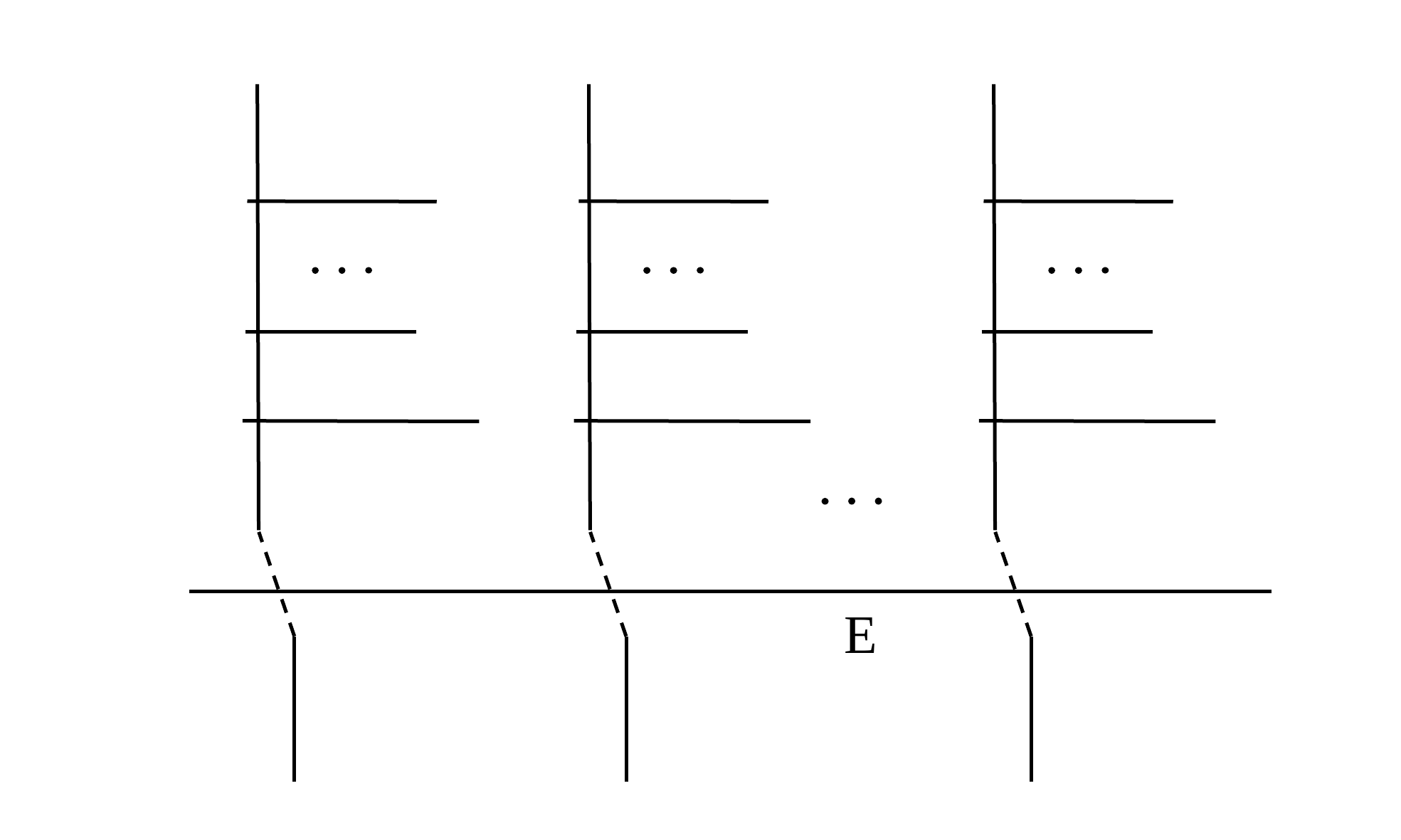}\caption{The divisor $D$ on $X$. Black lines are chains of rational curves. $(-1)$-curves are dashed. The big horizontal line is $E$.}  \label{fig:Xbarprim}\end{figure}

The divisor $D$ is a tree of smooth rational curves and the number of its components reflects the complexity of singularities of $\bar E$. Because of the minimality of the resolution, the exceptional divisor over each cusp, which is a part of $D$, contains at least two maximal twigs of $D$, hence $D$ has at least $2c$ maximal twigs (see Fig.\ \ref{fig:Xbarprim}). Based on the above theorem we prove the following result.

\bthm\label{thm:MAIN2} Let $\bar E\subseteq \PP^2$ be a complex rational cuspidal curve and let $(X,D)\to (\PP^2,\bar E)$ be the minimal log resolution of singularities. If any of the following conditions holds:\benum[(a)]

\item $\bar E$ has more than two cusps,

\item $D$ has more than six maximal twigs,

\item $\PP^2\setminus\bar E$ is not of log general type,
\eenum
then $\bar E$ satisfies the Coolidge-Nagata conjecture, i.e.\ there exists a Cremona transformation of $\PP^2$ which maps $\bar E$ onto a line. \ethm

We also prove the following important property of the birational minimalization morphism $\psi'\:(X,\frac{1}{2}D)\to (X_n',\frac{1}{2}D_n')$ resulting from the log MMP (see Sections \ref{ssec:logMMP} and \ref{sec:n<=1}).

\bthm\label{thm:MAIN3} Under the assumptions of Theorem \ref{thm:MAIN1} the minimalization morphism $\psi'\:(X,\frac{1}{2}D)\to (X_n',\frac{1}{2}D_n')$ contracts at most one irreducible curve not contained in $D$, i.e.\ $n\leq 1$. \ethm

\ssk

We now explain our strategy. Let $(X,D)\to (\PP^2,\bar E)$ be as above, with $\bar E\subseteq \PP^2$ not Cremona equivalent to a line. The case $\kappa(X\setminus D)\leq 1$ is done using structure theorems for open surfaces (\ref{prop:Fujita_boundary_excluded}), so we may assume $\kappa(X\setminus D)=2$. But now $X\setminus D$ is a $\Q$-acyclic surface of log general type, so it is well known that it does not contain lines (after removing the line we would get an affine surface of log general type with non-positive Euler characteristic, which is impossible by \ref{lem:BMY}). Then the birational morphism onto the (singular) minimal model of $(X,D)$ is well described, it contracts only the rays supported on $D$ and its subsequent images (hence $(X,D)$ is \emph{almost minimal}, see \cite[2.3.11, 2.4.3]{Miyan-OpenSurf}). Unfortunately the nefness of the log canonical divisor of this minimal model, even coupled with the log BMY inequality, gives rather weak consequences. Our key idea in this situation is simple. By the Kumar-Murthy criterion \ref{lem:cuspidal_of_gt_and_khalf}(iii) we have $2K_X+E\geq 0$, so in particular $\kappa(K_X+\frac{1}{2}D)\geq \kappa(K_X+\frac{1}{2}E)\geq 0$. This means that we can run the log MMP for $(X,\frac{1}{2}D)$ to obtain a minimal model $(X',\frac{1}{2}D')$ for which $K_{X'}+\frac{1}{2}D'$ is nef. The occurrence of $\frac{1}{2}$ changes (complicates) the situation entirely. Now $(X,\frac{1}{2}D)$ is not any more almost minimal, i.e.\ some curves not contained in $D$ are contracted, but the advantage is that the nefness of $K_{X'}+\frac{1}{2}D'$ is a much stronger property. In particular, our main inequality (\ref{prop:(Xi,Di)_properties_for_CN}) is obtained simply by intersecting the latter divisor with the push-forward of the effective divisor $K_X+\frac{1}{2}E$. In fact we gain a complete control over the shape and weights of $D'$, because by \cite[1.2(2)]{Palka-minimal_models} the possible pairs $(X_n',D_n')$ are well described, in particular there is a finite list of possible weighted dual graphs of $D_n'$ and in our situation we have additional restrictions \ref{cor:b2_bounds_etc}. However, since the obtained bounds do not immediately lead to a contradiction, we need to translate them to the level of $(X,D)$ and rule out the remaining configurations by referring to the geometry of cuspidal curves on $\PP^2$. Not to loose anything in translation, in Section \ref{sec:logMMP}, which is the technical core of the article, we introduce tools to carefully control the process of minimalization. We show that the minimalization morphism $\psi'\:(X,\frac{1}{2}D)\to (X',\frac{1}{2}D')=(X_n',\frac{1}{2}D'_n)$ is a composition of $n\geq 0$ contractions $\psi_i'\:(X_{i-1}',\frac{1}{2}D_{i-1}')\to (X_i',\frac{1}{2}D_i')$, $i=1,\ldots,n$, such that $X_i'\setminus D_i'$ is an open subset of $X_{i-1}'\setminus D_{i-1}'$ with the complement isomorphic to $\C^*$. It turns out that as long as $D$ is not very small (has more than six maximal twigs) we are able to deal with the situation completely. Theorem \ref{thm:MAIN3} rules out all cases in which $n>1$. Thus the essential difficulty to overcome to prove the Coolidge-Nagata conjecture in the remaining cases of uni- and bi-cuspidal curves is the situation when $n=1$ and $p_2(\PP^2,\bar E)=3$ (see \ref{cor:p2=/=4}), or equivalently, $(K_X+D)^2=1$. We will address this problem in a forthcoming paper.

\med At a conference in Montreal in September 2012, \emph{The Topology of Algebraic Varieties}, Mariusz Koras, who is independently working on the Coolidge-Nagata conjecture, has shown a proof that $(K_X+D)^2\neq -1$ using different methods and announced that he proved the conjecture for unicuspidal curves. At the same conference, the author has shown how to extend the methods from \cite{Palka_CN-arxiv} to prove it for cuspidal curves with more than three cusps. The ideas from \cite{Palka_CN-arxiv} are used in Section \ref{sec:CN} after the corollaries from the log MMP are established.

\ssk\textsl{\textsf{Acknowledgements.}} The author would like to thank Mikhail Zaidenberg for turning his attention to cuspidal curves and discussing results in the literature.

\tableofcontents

\section{Preliminaries}\label{sec:preliminaries}

We recall basic definitions from the theory of non-complete surfaces (see \cite{Fujita-on_the_topology} and \cite{Miyan-OpenSurf}). We work over the field of complex numbers.

\subsection{Log surfaces and divisors}\label{ssec:open surfaces}
Given two $\Q$-divisors $T$, $T'$ we say that $T'$ is a \emph{subdivisor of $T$}, we write $T'\leq T$, if $T-T'$ is effective.  Let $T$ be a reduced divisor on a smooth projective surface $X$. If $R$ is a reduced subdivisor of $T$ we define $\beta_T(R)=R\cdot (T-R)$ and we call it a \emph{branching number of $R$ in $T$}. If $R$ is irreducible and nonzero we say that $R$ is a \emph{tip} or a \emph{branching component} if $\beta_T(R)\leq 1$ or $\beta_T(R)\geq 3$ respectively. We say that $T$ is an snc-divisor if its irreducible components are smooth and intersect transversally, at most two in one point. We call $T$ a \emph{chain} if it is a connected snc-divisor whose dual graph is linear. A subdivisor of $T$ is a \emph{twig} of $T$ if it is a chain which contains no branching components of $T$ and contains a tip of $T$. We say that $T$ is a \emph{(rational) tree} if it is an snc-divisor with connected support (whose all components are rational and) such that the dual graph of $T$ contains no loops. A fork is a tree with exactly one branching component. The \emph{arithmetic genus of $T$} is $p_a(T)=\frac{1}{2}T\cdot (K+T)+1,$ where $K$ is the canonical divisor (class) on $X$. For a rational tree $p_a(T)=0$. We denote the Iitaka-Kodaira dimension of the divisor $T$ on $X$ by $\kappa(T)$ and the Picard rank, i.e.\ the rank of the Neron-Severi group of $X$ by $\rho(X)$. If $T=T_1+\ldots+T_k$ is a decomposition of an ordered rational chain into irreducible components, such that $T_i\cdot T_{i+1}=1$ for $i<k$, then we write $T=[-T_1^2,-T_2^2,\ldots,-T_k^2].$ By $(m)_p$ we mean a sequence $(m,m,\ldots,m)$ of length $p$. By a curve we mean a one-dimensional irreducible and reduced variety. An \emph{$(n)$-curve} is a smooth rational curve with self-intersection $n$. A $(-1)$-curve which is a component of $T$ as above is called \emph{superfluous} if it intersects at most two other components of $T$, each at most once and transversally. We define the \emph{discriminant of $T$} as $d(T)=\det(-Q(T))$, where $Q(T)$ is the intersection matrix of $T$. We put $d(0)=1$. By $\#T$ we denote the number of irreducible components of $T$.

If $T'$ is a rational twig without superfluous $(-1)$-curves and with a negative definite intersection matrix, or more generally, if it is a rational ordered chain (i.e.\ a chain with a choice of a tip) with these properties, we put $$\ind (T')=\frac{d(T'-\tip(T'))}{d(T')} \text{\ \ and\ \ } \delta(T')=\frac{1}{d(T')}.$$ The former number is usually called the \emph{inductance} or \emph{capacity} of $T'$.

Assume now that $T$ is a connected snc-divisor with rational components and no superfluous $(-1)$-curves. Assume that its intersection matrix is not negative definite (this is the case when $T$ is an snc-minimal boundary of an affine surface) and that the intersection matrices of all its maximal twigs are negative definite (this is the case if $\kappa(K+T)\geq 0$). Let $T_i$, $i=1,\ldots,t$, be the maximal twigs of $T$. We put $$\ind (T)=\sum_{i=1}^t \ind (T_i)\text{\ \ and \ } \delta(T)=\sum_{i=1}^t\delta(T_i) .$$

Assume additionally that $\kappa(K+T)\geq 0$. We have the Fujita-Zariski decomposition $K+T=(K+T)^++(K+T)^-,$ where $(K+T)^+$ is numerically effective and $(K+T)^-$ is effective, either empty or having a negative definite intersection matrix. Moreover, $(K+T)^+\cdot B=0$ for any curve $B$ contained in $\Supp (K+T)^-$. If $T_i$ is a (not necessarily maximal) twig of $T$ then we define $\Bk_T T_i$, the \emph{bark} of $T_i$ with respect to $T$, as the unique $\Q$-divisor supported on $\Supp T_i$, such that $$\Bk_T T_i\cdot R=\beta_T(R)-2$$ for every component $R$ of $T_i$, equivalently that $\Bk_T T_i\cdot R$ equals $-1$ if $R$ is the tip of $T$ contained in $T_i$ and is zero otherwise. Then we define $\Bk T$, the \emph{bark of $T$}, as the sum of $\Bk_T T_i$'s taken over all maximal twigs of $T$. Note that if $\alpha\:(X,T)\to (X',\alpha_*T)$ is the contraction of maximal twigs of $T$ then $$\alpha^*(K_{X'}+\alpha_*T)=K_X+T-\Bk T.$$ A $(-2)$-twig is a twig whose all irreducible components are $(-2)$-curves. A maximal $(-2)$-twig is a $(-2)$-twig which is not a proper subdivisor of another $(-2)$-twig. The following lemma summarizes what we need to know about barks.

\blem\label{lem:Bk} Let $T_i$ be a maximal twig of a rational tree $T$ as above and let $T_0$ be a component of $T_i$. Denote the coefficient of $T_0$ in $\Bk T$ by $t_0$. Then \benum[(i)]

\item $0<t_0<1$.

\item If $T_0$ is the component meeting $T-T_i$ then $t_0=\delta(T_i)$.

\item $(\Bk T)^2=-\ind (T)$.

\item If there is no $(-1)$-curve $A$ on $X$, for which $T\cdot A\leq 1$, then $(K+T)^-=\Bk T$.

\eenum \elem

\begin{proof} Write $T_i=T_{i,1}+T_{i,2}+\ldots+T_{i,k_i}$, where $T_{i,j}$ are irreducible and $T_{i,j}\cdot T_{i,j+1}=1$ for $j<k_i$. Then by \cite[2.3.3.4]{Miyan-OpenSurf} the coefficient of $T_{i,j}$ in $\Bk T$ equals $d(T_{i,j+1}+\ldots+T_{i,k_i})/d(T_i)$. This gives (i), (ii) and (iii). Part (iv) follows from \cite[6.20]{Fujita-on_the_topology}.
\end{proof}

We will often use the logarithmic version of the Bogomolov-Miyaoka-Yau inequality proved originally by Kobayashi, Nakamura and Sakai. We use the version due to Langer \cite{Langer}. The following formulation is taken from \cite[2.5]{Palka-exceptional}. Recall that a connected reduced divisor $D$ is of \emph{quotient type} if it can be contracted to a quotient singularity. We denote by $\Gamma(D)$ the local fundamental group of the corresponding singular point. For the definition of almost minimality see Section \ref{sec:main result}.

\begin{lemma}[The log BMY inequality]\label{lem:BMY}
Let $(X,D)$ be a pair consisting of a smooth projective surface $X$ and a reduced snc-divisor $D$. Let $D_1,\dots, D_k$ be the connected components of $D$ which are of quotient type. If $(X,D)$ is almost minimal and $\kappa(X\setminus D)\geq 0$ then \[\frac{1}{3}(K_{X}+D-\Bk D)^2\leq \chi(X\setminus D)+\sum_{i=1}^k\frac{1}{|\Gamma(D_i)|}.\] \end{lemma}

For a pair $(X,D)$ as in the lemma let $\sigma\:(X',D')\to (X,D)$, where $D'=\sigma_*^{-1}D+\Exc \sigma$, be a blowup with a center on $D$. Write $\sigma^*D=D'+\mu\Exc \sigma$. We say that $\sigma$ is \emph{inner (outer)} for $D$ if $\mu=1$ ($\mu=0$ respectively). We may equivalently ask that the center of $\sigma$ belongs to exactly two (one) components of $D$. If $C$ is an irreducible curve on $X'$ then we say that $\sigma$ \emph{touches $C$} if and only if $\Exc \sigma\cdot C\neq 0$. If $\Exc \sigma\cdot C=1$ we say that $\sigma$ \emph{touches $C$ once}. In a more general situation, when $\sigma \:X'\to X$ is a birational morphism from a smooth projective surface we say that $\sigma$ \emph{touches $C$ $n$ times} if in the decomposition of $\sigma$ into blowups each time the respective image of $C$ is touched at most once, and is touched $n$ times in total. If we say that two effective divisors on $X$ \emph{meet $n$ times} we simply mean that their intersection number equals $n$.

By a \emph{log surface} we mean a pair $(Y,B)$ consisting of a projective normal surface $Y$ together with an effective $\Q$-divisor, which can be written as $B=\sum b_iB_i$, where $B_i$ are distinct irreducible components and $0<b_i\leq 1$. It is \emph{smooth} if $X$ is smooth and $\sum B_i$ is an snc-divisor.

\bdfn\label{def:resolution} Let $(Y,B)$ be a log surface and let $\pi\:(X,D)\to (Y,B)$ be a proper birational morphism from a log surface such that $X$ is smooth and $D=\pi_*^{-1}B+\Exc \pi$. We say that $\pi$ is a \emph{weak (embedded) resolution of singularities} if $\pi_*^{-1}B$ is an snc-divisor. It is a \emph{log resolution} if $D$ is an snc-divisor (equivalently, $(X,D)$ is a smooth log surface). A log (resp.\ weak) resolution is a \emph{minimal} log (resp.\ weak) resolution if it does not dominate any other log (resp.\ weak) resolution.\edfn

\subsection{Rational cuspidal curves}\label{ssec:cuspidal_curves}

For a rational curve $\bar E\subseteq \PP^2$ we put $\kappa_{1/2}(\PP^2,\bar E)=\kappa(K_X+\frac{1}{2}D)$ and $p_2(\PP^2,\bar E)=h^0(2K_X+D)$, where $(X,D)\to (\PP^2,\bar E)$ is the minimal log resolution.

\blem\label{lem:cuspidal_of_gt_and_khalf} Let $\bar E\subseteq \PP^2$ be a rational cuspidal curve and let $\pi\:(X,D) \to (\PP^2,\bar E)$ be any weak resolution of singularities. Put $E=\pi_*^{-1}\bar E$. Then: \benum[(i)]

\item $\PP^2\setminus\bar E$ is $\Q$-acyclic.

\item $h^0(2K_X+D)=p_2(\PP^2,\bar E)$.

\item $\bar E\subseteq \PP^2$ is Cremona equivalent to a line if and only if $h^0(2K_X+E)=0$. In particular, if $\bar E\subseteq \PP^2$ is not Cremona equivalent to a line then $E^2\leq -4$.

\item If $\bar E\subseteq \PP^2$ is not Cremona equivalent to a line and $f$ is a fiber of any $\PP^1$-fibration of $X$ then $f\cdot E\geq 4$.

\item If $\PP^2\setminus \bar E$ is of log general type then it contains no topologically contractible curves. In case $\pi$ is the minimal log resolution we have $(K_X+D)^-=\Bk D$.

\item Let $\mu(p)$ denote the multiplicity of a point $p\in\bar E$. Then
$$\deg \bar E<3\max \{\mu(p):p\in\Sing \bar E\}.$$

\eenum \elem

\begin{proof} (i) The $\Q$-acyclicity follows from the Lefschetz duality.

(ii)  Let $\sigma\:X'\to X$ be a blowup and let $D'=\sigma_*^{-1}D+\Exc \sigma$. Let $\mu$ be the number of components of $D$ passing through the center of $\sigma$. The divisor $D-E$ is snc and $E$ is smooth, so $\mu\leq 3$. Clearly, $\sigma_*$ embeds the linear system of $2K_{X'}+D'$ into the linear system of $\sigma_*(2K_{X'}+D')=2K_X+D.$ Now if $2K_X+D\sim U$ then $2K_{X'}+D'\sim \sigma^*U+(3-\mu)\Exc \sigma\geq \sigma^*U,$ so $\sigma_*$ is surjective. It follows that $h^0(2K_X+D)$ does not depend on the choice of a weak resolution.

(iii) The first statement is a consequence of \cite[2.4, 2.6]{Kumar-Murthy} and holds for any rational curve in $\PP^2$. If $2K+E\geq 0$ and $E^2\geq -3$ then $E\cdot (2K+E)=2E\cdot (K+E)-E^2=-4-E^2<0$, so $E$ is in the fixed part of $2K+E$, which contradicts the rationality of $X$.

(iv) Let $f$ be a smooth fiber of a $\PP^1$-fibration of $X$. By (ii) and (iii) $h^0(2K_X+E)>0$, so $0\leq f\cdot (2K_X+E)=-4+f\cdot E$.

(v) By \cite{MiTs-lines_on_qhp} $\PP^2\setminus\bar E$ contains no topologically contractible curves. If $\pi$ is the minimal log resolution then by \ref{lem:Bk}(iv) $(K_X+D)^-=\Bk D$.

(vi) This is the inequality of Matsuoka-Sakai \cite{MaSa-cusp}.
\end{proof}

The following result was shown in \cite[2.6]{Palka-minimal_models}. We recall the proof for completeness.

\bprop\label{prop:Fujita_boundary_excluded} If $\PP^2\setminus \bar E$ is not of log general type then it is $\C^1$- or $\C^*$-fibered and $\kappa_{1/2}(\PP^2,\bar E)=-\8$. In particular, $p_2(\PP^2,\bar E)=0$ and $\bar E$ satisfies the Coolidge-Nagata conjecture. \eprop

\begin{proof} Let $S=\PP^2\setminus\bar E$. If $\kappa(S)=-\8$ or $1$ then $S$ is $\C^1$- or $\C^*$-fibered by structure theorems for smooth affine surfaces (see \cite[3.1.3.2, 3.1.7.1]{Miyan-OpenSurf}), so we may assume $\kappa(S)=0$. Let $(X,D)$ be the minimal log resolution of $(\PP^2,\bar E)$.

Assume first that $(X,D)$ is not almost minimal. Since $D$ is connected, by \cite[6.20]{Fujita-on_the_topology} $S$ contains a curve $\ell$ isomorphic to $\C^1$ such that $\kappa(S\setminus \ell)=\kappa(S)=0$. Since $\Pic (S)$ is torsion, there is a rational map $f\:S\map \C^1$ such that $(f)=m\ell$ for some positive integer $m$. Because $S$ is affine we may assume that $f$ is regular on $S$, so we get a morphism $f\:S\to \C^1$. If $f'\:S\to B$ is its Stein factorization then $\kappa(B-f'(\ell)) \geq\kappa(\C^*)=0$ and $0=\kappa(S\setminus\ell)\geq \kappa (F)+ \kappa(B-f'(\ell))$ for a general fiber $F$ over $B-f'(\ell)$ by the Kawamata addition theorem \cite{Kawamata}. Since $S\setminus \ell$ is not $\C^1$-fibered, we get $\kappa(F)=0$, so $f$ is a $\C^*$-fibration and we are done.

Suppose now that $(X,D)$ is almost minimal but $S$ is not $\C^*$-fibered. Because $\kappa(S)=0$, by \cite[8.64]{Fujita-on_the_topology} (cf. \cite[3.4.4.2]{Miyan-OpenSurf}) $\PP^2\setminus \bar E$ is one of the three Fujita surfaces $Y\{a,b,c\}$. Then $D$ is a fork whose maximal twigs are $(-2)$-chains. But then $c=1$ and there should exist a tip of $D$ (namely $E$) for which $D-E$ is negative definite. This is not so for the surfaces $Y\{a,b,c\}$, because the branching components have self-intersections $1,0,-1$ respectively; a contradiction.

Thus $S$ is $\C^1$- or $\C^*$-fibered. The fibration extends to a $\PP^1$-fibration of some weak resolution $(X,D)\to (\PP^2,\bar E)$, so that $f\cdot D\leq 2$. But then $\kappa_{1/2}(\PP^2,\bar E)=\kappa(K_X+\frac{1}{2}D)=-\8$, because $f\cdot(K_X+\frac{1}{2}D)<0$. It follows that $p_2(\PP^2,\bar E)=h^0(2K_X+E)=0$, so we are done due to \ref{lem:cuspidal_of_gt_and_khalf}(iii).
\end{proof}

\brem\label{rem:deg<=5_rectifiable} Let $\bar E\subseteq\PP^2$ be a rational curve of degree $d$ which is not Cremona equivalent to a line. The criterion of Kumar-Murthy stated in \ref{lem:cuspidal_of_gt_and_khalf}(iii) for cuspidal curves works for all rational curves, so $\kappa(K_X+\frac{1}{2}E)\geq 0$ and therefore $\kappa(K_X+\frac{1}{2}D)\geq 0$. We have the equality $\pi_*(2K_X+E)=2K_{\PP^2}+\bar E$, which implies that $d\geq 6$. For a partial converse note that if $\bar E\subseteq \PP^2$ is a general rational curve of degree $d$ then its singularities are ordinary double points (nodes), so we compute easily $2K_X+E\sim (d-6)H$, where $H$ is the pullback of a line on $\PP^2$. Thus such $\bar E$ is not Cremona equivalent to a line for $d\geq 6$. A general rational sextic with ten nodes is an example of lowest degree. On the other hand, if the singularities of $\bar E$ are more complicated than nodes (like for cuspidal curves) determining whether $\bar E$ is Cremona equivalent to a line is a more subtle issue. \erem

\section{Characteristic pairs}

To describe the geometry of exceptional divisors for a resolution of singularities we rely on the Hamburger-Noether pairs. We follow \cite[1.2]{Koras-ab_moh} (for more details see \cite{Russell2}). Some authors prefer to work with Puiseux expansions and Puiseux pairs, which can be used to do analysis analogous to ours (see \cite{OrZa_cusp}). We note that one of the advantages of Hamburger-Noether pairs, aside of their explicit geometric interpretation, is that they can be used when working over fields of positive characteristic.

As an input data take an analytically irreducible germ of a singular curve $(\chi,q)$ on a smooth surface and a (germ of a) curve $C$ passing through $q$, smooth at $q$. Put $(C_1,\chi_1,q_1)=(C,\chi,q)$, $c_1=(C_1\cdot \chi_1)_{q_1}$, where $(\ \cdot\ )_{q_1}$ denotes the local intersection index at $q_1$, and choose a local coordinate $y_1$ at $q_1$ in such a way that $Y_1=\{y_1=0\}$ is transversal to $C_1$ at $q_1$ and $p_1=(Y_1\cdot\chi_1)_{q_1}$ is not bigger than $c_1$. Blow up over $q_1$ until the proper transform $\chi_2$ of $\chi_1$ intersects the reduced total transform of $C_1+Y_1$ not in a node. Let $q_2$ be the point of intersection and let $C_2$ be the last exceptional curve. Put $c_2=(C_2\cdot \chi_2)_{q_2}$. We repeat this procedure and we define successively $(\chi_i,q_i)$ and $C_i$ until $\chi_{h+1}$ is smooth for some $h\geq 1$. This defines a sequence $\binom{c_1}{p_1}, \binom{c_2}{p_2}, \ldots, \binom{c_h}{p_h}$, depending on the choice of $C$. It follows from the definition that $c_1\geq p_1$ and that $p_1$ is the first (and maximal) number in the sequence of multiplicities of $q\in \chi$. Note that the total exceptional divisor contains a unique $(-1)$-curve.

Because of the forced condition $p_i\leq c_i$ the sequence is usually longer than the sequence of Puiseux pairs. Although it is defined for (and depends on) any initial curve $C$, in this article we will choose for $C$ a smooth germ maximally tangent to $\chi$ (note that because $\chi$ is singular its intersection with smooth germs passing through $q$ is bounded from above). For this choice of $C$ we refer to the sequence $\binom{c_1}{p_1}, \binom{c_2}{p_2}, \ldots, \binom{c_h}{p_h}$ as the sequence of \emph{Hamburger-Noether pairs} (or \emph{characteristic pairs} for short) \emph{of the resolution of $(\chi_1,q_1)$} (as it follows from the definition we refer here to the minimal log resolution). It is convenient to extend the definition to the case when $(\chi_1,q_1)$ is smooth by defining its sequence of characteristic pairs to be $\binom{1}{0}$. By $\binom{c}{p}_k$ we mean a sequence of pairs $\binom{c}{p},\ldots,\binom{c}{p}$ of length $k$. Let $(\mu_i)_{i\in I_j}$ be the non-increasing sequence of multiplicities of successive centers for the sequence of blowups as above leading from $\chi_j$ to $\chi_{j+1}$. The sequence $(\mu_i)_{i\in I_1},\ldots,(\mu_i)_{i\in I_h}$ is the \emph{multiplicity sequence of the singularity $(\chi,q)$}. Note that the composition of blowups corresponding to multiplicities bigger than $1$ is the minimal weak resolution of singularities.

Let now $\pi\:X\to X'$ be a proper birational morphism of smooth surfaces, such that the exceptional divisor $Q=\Exc \pi$ contains a unique $(-1)$-curve $U$. If $U$ is not a tip of $Q$ then $\pi$ is a minimal log resolution of some (germ of a) singular curve $\chi$ on $X'$ as above and we define the \emph{sequence of characteristic pairs of $Q$} to be the one of $\chi$. In case $U$ is a tip of $Q$ let $(X,Q)\to (Y,Q')$ be a composition of a minimal number of contractions, say $m$, of $(-1)$-curves in $Q$ and its successive images, such that $Q'$ contains no $(-1)$-tip. If $(\binom{c_i}{p_i})_{i\leq h}$ is the sequence of characteristic pairs for $Q'$ then the sequence of characteristic pairs of $Q$ is by definition $(\binom{c_i}{p_i})_{i\leq h},\binom{1}{1}_m$.

\blem\label{1lem:sum of mu and squares} Assume that the sequence of blowups $(\sigma_j)_{j\in I_i}$, leading from $(\chi_i,q_i)$ to $(\chi_{i+1},q_{i+1})$ is described as above by the characteristic pair $\binom{c_i}{p_i}$. Let $\mu_j$ be the multiplicity of the center of $\sigma_j$ as a point on the proper transform on $\chi_i$. Then we have:\benum[(i)]

\item $c_{i+1}=gcd(c_i,p_i)$,

\item  $\ds \sum_{j\in I_i} \mu_j=c_i+p_i-gcd(c_i,p_i)$,

\item $\ds \sum_{j\in I_i}\mu_j^2=c_ip_i$. \eenum \elem

\begin{proof}The formulas hold in case $c_i=p_i$. If $c_i>p_i$ then perform the first blowup and note that the remaining part of the sequence $(\sigma_j)_{j\in I_i}$ is described by $\binom{c_i-p_i}{p_i}$ in case $c_i-p_i\geq p_i$ or by $\binom{p_i}{c_i-p_i}$ otherwise. The multiplicity of the first center is $p$. Now the result follows by induction on $\max(c_i,p_i)$.
\end{proof}

Let $\bar E\subseteq \PP^2$ be a rational cuspidal curve with cusps $q_1,\ldots,q_c$ and let $\rho\:(Y,B)\to (\PP^2,\bar E)$, with $B=\rho_*^{-1}\bar E+\Exc \rho$, be a weak resolution of singularities, such that each exceptional divisor $Q_i=\rho^{-1}(q_i)_{red}$ contains a unique $(-1)$-curve. Put $E_Y=\rho_*^{-1}\bar E$. For $i=1,\ldots,c$ let $\bar \mu_i=(\mu_{i,1},\ldots,\mu_{i,k_i})$ be the multiplicity sequence of $q_i\in \bar E$ and its infinitely near points for this resolution. Denote the sequence of characteristic pairs of the divisor $Q_i$ by $(\binom{c_{i,j}}{p_{i,j}})_{j\leq h_i}$. Note that the multiplicity of $q_i\in \bar E$ is $\mu(q_i)=\mu_{i,1}=p_{i,1}$.

\bcor\label{cor:sum of m_i, m_i^2 using pairs} Assume $E_Y$ meets $Q_i$ not in a node. Then: \benum[(i)]

\item $\ds \sum_{j=1}^{k_i} \mu_{i,j}=(Q_i\cdot E_Y) (c_{i,1}+\sum_{j=1}^{h_i}p_{i,j}-1)$,

\item $\ds \sum_{j=1}^{k_i} \mu_{i,j}^2=(Q_i\cdot E_Y)^2 \sum_{j=1}^{h_i}c_{i,j}p_{i,j}.$

\eenum \ecor

\bcor\label{cor:I and II} Let the situation be as above. Put $\gamma_Y=-E_Y^2$, $d=\deg \bar E$, $\rho_i=Q_i\cdot E_Y$, $M(q_i)=c_{i,1}+\ds \sum_{j=1}^{h_i}p_{i,j}-1$ and $I(q_i)=\ds \sum_{j=1}^{h_i}c_{i,j}p_{i,j}$. Then \benum[(i)]

\item $\gamma_Y-2+3d=\ds \sum_i \rho_i M(q_i),$

\item $\gamma_Y+d^2=\ds \sum_i \rho_i^2 I(q_i),$

\item $(d-1)(d-2)=\ds \sum_i \rho_i(\rho_i I(q_i)-M(q_i)).$
\eenum \ecor

Note that (iii) is the genus-degree formula written in terms of characteristic pairs.

\begin{proof} Let $C\subseteq X$ be an irreducible curve on a smooth projective surface. Let $p\in C$ be a singular point of $C$ having multiplicity $m$ and let $\sigma\: X'\to X$ be a blowup at $p$. Denote the exceptional curve by $L$ and the proper transform of $C$ on $X'$ by $C'$. Then $$K_X\cdot C=\sigma^*K_X\cdot C'=(K_{X'}-L)\cdot C'=K_X'\cdot C'-m$$ and $$C^2=\sigma^*C\cdot C'=(C'+mL)\cdot C'=C'^2+m^2.$$ This implies that if in the lemma $\rho_i=1$ then the sum of all multiplicities $\mu_{i,j}$ equals $K_Y\cdot E_Y-K_{\PP^2}\cdot \bar E=\gamma_Y-2+3d$ and the sum of their squares equals $\bar E^2-E_Y^2=d^2+\gamma_Y$. The general case follows by the linearity of multiplicity with respect to addition of germs. This proves (i) and (ii); (iii) is their consequence.
\end{proof}

\bex If $Q_i=[2,1,3]$ then $h_i=1$, $\binom{c_{i,1}}{p_{i,1}}=\binom{3}{2}$, $M(q_i)=4$ and $I(q_i)=6$. If $Q_i=[5,2,1,3,2,2,3]$ then $h_i=1$, $\binom{c_{i,1}}{p_{i,1}}=\binom{16}{9}$, $M(q_i)=24$ and $I(q_i)=144$.  The multiplicity sequence for a cusp with $Q_i=[2,1,3]$ is $(2,1,1)$. For $Q_i=[3,1,2,3]$ it is $(3,2,1,1)$, for $Q_i=[n,1,(2)_{n-2},3]$ it is $(n,n-1,(1)_{n-1})$ and for $Q_i=[5,2,1,3,2,2,3]$ it is $(9,7,(2)_3,(1)_2)$.\eex

\blem\label{lem:ind_using_pairs} Let $Q$ be a rational chain created by the characteristic pair $\binom{c}{p}$ and let $L$ be the unique $(-1)$-curve of $Q$. Write $Q-L=A+B$, where $A$ and $B$ are disjoint and connected with $d(A)\geq d(B)$. Put $B=0$ if $p=1$. Let $A_L, B_L$ be the tips of $Q$ which are components of $A$ and $B$ respectively. Write $c=q\cdot p+r$ for some integers $q>0$ and $0\leq r<p$. Then $$d(A)=c,\ d(B)=p,\ d(A-A_L)=c-p \text{\ \ and\ \ } d(B-B_L)=p-r.$$ Moreover, if we treat $A$ and $B$ as twigs of $Q$ then
$$\ind (A)+\ind (B)=\frac{c-p}{c}+\frac{p-r}{p}>1-\frac{1}{q+1}\geq \frac{1}{2}.$$ \elem

\begin{proof} For an  (ordered) rational chain $T=T_1+\ldots+T_n$ we put $d'(T)=d(T-T_1)$ if $n>0$ and $d'(T)=0$ if $n=0$. As a consequence of elementary properties of determinants we have a recurrence formula computing the discriminant $$d(T)=(-T_1^2)\cdot d(T-T_1)-d'(T-T_1).$$ We assume that $A$ and $B$ are ordered so that the tips of $Q$ are their first components. We prove the four formulas for discriminants by induction on the length of $Q$. Let $\psi\: X\to X'$ be a birational morphism for which $\Exc \psi=Q$ and let $\sigma$ be the blowup with the center $\psi(Q)$. Decompose $\psi$ as $\psi=\sigma \circ \psi'$ and let $Q'=\Exc \psi'$. The proper transform $U$ of $\Exc \sigma$ is a component of $A$ which is a tip of $Q$. We have $Q-U=Q'$. If $q>1$ then $U$ is a $(-2)$-curve and $Q'$ is created by the characteristic pair $\binom{c-p}{p}$. Then by the inductive assumption we have $d(A-A_L)=c-p$, $d(B)=p$ and $d(B-B_L)=p-r$. We compute also $$d(A)=2d(A-U)-d'(A-U)=2(c-p)-(c-2p)=c,$$ so we are done. Assume $q=1$. Let $q'$ and $r'$ be the quotient and the remainder of dividing $p$ by $c-p$. Then $U$ is a $-(q'+2)$-curve and $Q'$ is created by the characteristic pair $\binom{p}{c-p}$. By the inductive assumption $d(B)=p$, $d(B-B_L)=p-(c-p)=p-r$ and $d(A-A_L)=c-p$. We compute $$d(A)=(q'+2)(c-p)-(c-p-r')=q'(c-p)+r'+(c-p)=c.$$

We have $$\frac{c-p}{c}+\frac{p-r}{p}>1-\frac{1}{q+1} \iff \frac{p-r}{p}+\frac{1}{q+1}>\frac{p}{c} \iff \frac{(p-r)(q+1)+p}{p}>\frac{c-r+p}{c} \iff$$ $$\iff \frac{(p-r)(q+1)}{p}>\frac{p-r}{c} \iff \frac{c}{p}>\frac{1}{q+1}.$$ But $\frac{c}{p}>1>\frac{1}{q+1}$, so we are done.
\end{proof}

\blem\label{lem:Q_with_small_KQ} Let $Q$ be a rational chain which contracts to a smooth point and which contains a unique $(-1)$-curve. Let $k$ denote a non-negative integer. \benum[(i)]

\item If $K\cdot Q<0$ then $Q=[(2)_k,1]$.

\item If $K\cdot Q=0$ then $Q=[(2)_k,3,1,2]$.

\item If $K\cdot Q=1$ then $Q=[(2)_k,4,1,2,2]$ or $Q=[(2)_k,3,2,1,3]$.

\item If $K\cdot Q=2$ then $Q=[(2)_k,5,1,2,2,2]$ or $Q=[(2)_k,4,2,1,3,2]$ or $Q=[(2)_k,3,3,1,2,3]$ or $Q=[(2)_k,3,2,2,1,4]$.
\eenum
\elem

\begin{proof} If $\#Q=1$ then $Q=[1]$, so $K\cdot Q=-1$. Let $\sigma$ be a blowup with a center on $Q$ and let $Q'$ be its reduced total transform. If the blowup is outer then $K\cdot Q'=K\cdot Q$. It the blowup is inner then $K\cdot Q'=K\cdot Q+1$. Therefore, every chain with $K\cdot Q=\alpha\geq 0$ can be obtained from some chain with $K\cdot Q=\alpha-1$ by blowing up once. Clearly, for inner blowups there are two possible choices of the center.
\end{proof}

More generally, we have the following description.

\blem\label{lem:chains_with_d=1} Every rational chain $Q$ which contracts to a smooth point and which contains a unique $(-1)$-curve is of the following type: $$[m_{2k}+3,(2)_{m_{2k-1}},\ldots,m_2+3,(2)_{m_1},1,m_1+2,(2)_{m_2},m_3+3,\ldots,m_{2k-1}+3,(2)_{m_{2k}},3,(2)_x]$$ for some integers $k\geq 1$, $x\geq -1$ and $m_1,\ldots,m_{2k}\geq 0$ or $$[(2)_{m_{2k}},m_{2k-1}+3,\ldots,m_3+3,(2)_{m_2},m_1+2,1,(2)_{m_1},m_2+3,\ldots,(2)_{m_{2k-1}},m_{2k}+3,(2)_x]$$ for some integers $k,x\geq 0$ and $m_1,\ldots,m_{2k}\geq 0$. By definition $[3,(2)_{-1}]=\emptyset$.
\elem

\begin{proof} Induction with respect to the length of the chain.
\end{proof}

\section{Minimal models}\label{sec:logMMP}

From now on we assume that $\bar E\subseteq \PP^2$ is a rational cuspidal curve, $\pi_0\:(X_0,D_0)\to (\PP^2,\bar E)$ is the minimal weak resolution of singularities and $\pi\:(X,D)\to(\PP^2,\bar E)$ is the minimal log resolution. We assume also that $\bar E\subseteq \PP^2$ violates the Coolidge-Nagata conjecture, i.e.\ it is not Cremona equivalent to a line. By \ref{prop:Fujita_boundary_excluded} and \ref{lem:cuspidal_of_gt_and_khalf}(iii) $\PP^2\setminus \bar E$ is a surface of log general type and $2K+E\geq 0$, where $K$ denotes the canonical divisor on $X$. In particular, $p_2(\PP^2,\bar E)=h^0(2K+D)\geq 1$. Clearly, we have a factorization $\psi_0\:(X,D)\to (X_0,D_0)$. Our method to analyze $(\PP^2,\bar E)$ is to run the logarithmic Minimal Model Program for $(X_0,\frac{1}{2}D_0)$, which is described by $\alpha_0$ and the lower row in the diagram below. But since the surfaces $Y_i$ are singular it is worth keeping track of lifts of extremal contractions of rays not contained in the boundary on the level of minimal weak and log resolutions dominating successive surfaces produced by the program (the upper and middle rows of the diagram). The analysis is an analog of the theory of peeling \cite[2.3]{Miyan-OpenSurf} (which works for reduced boundaries).

\subsection{Review of the log MMP for $(X_0,\frac{1}{2}D_0)$}\label{ssec:logMMP} The outcome of the construction in \cite[\S3]{Palka-minimal_models} is the following commuting diagram:

$$ \xymatrix{
{} & {(X,D)}\ar@{>}[r]^-{\psi_1'}\ar@{>}[d]^{\psi_0}\ar@{>}[ld]_-{\ \pi} &{(X_1',D_1')}\ar@{>}[r]^-{\psi_2'\ }\ar@{>}[d]^{\varphi_1} &{\ldots}\ar@{>}[r]^-{\psi_n'\ } &{(X_n',D_n')}\ar@{>}[d]^{\varphi_n}\\
{(\PP^2,\bar E)} & {(X_0,D_0)}\ar@{>}[r]^-{\psi_1}\ar@{>}[d]^{\alpha_0}\ar@{>}[l]^-{\ \pi_0} &{(X_1,D_1)}\ar@{>}[r]^-{\psi_2\ }\ar@{>}[d]^{\alpha_1} &{\ldots}\ar@{>}[r]^-{\psi_n\ } &{(X_n,D_n)}\ar@{>}[d]^{\alpha_n}
\\ {} & {(Y_0,D_{Y_0})}\ar@{>}[r]^-{\psi_1''} & {(Y_1,D_{Y_1})}\ar@{>}[r]^-{\psi_2''\ } & {\ldots}\ar@{>}[r]^-{\psi_n''\ } & {(Y_n,D_{Y_n})}}$$

\med Let us recall the necessary definitions and results. Write $K_i$ (resp.\ $K_i'$) for the canonical divisor on $X_i$ (resp.\ on $X_i'$). Let $\varphi_i$, $i\geq 0$ be the minimal log resolution of $(X_i,D_i)$. We have $\varphi_0=\psi_0$. The morphism $\alpha_i\:(X_i,D_i)\to (Y_i,D_{Y_i})$ is the composition of successive extremal birational contractions of curves with support in $D_i$ which are negative with respect to $K_i+\frac{1}{2}D_i$ (and its successive pushforwards). The one-step description of $\alpha_i\:X_i\to Y_i$ is that it is the contraction of $\Delta_i+\Upsilon_i$, where $\Delta_i$ is the sum of all maximal $(-2)$-twigs of $D_i$ and $\Upsilon_i$ is the sum of $(-1)$-curves $L$ in $D_i$, for which either $\beta_{D_i}(L)=3$ and $L\cdot \Delta_i=1$ or $\beta_{D_i}(L)=2$ and $L$ meets exactly one component of $D_i$. Write $\Delta_i=\Delta^+_i+\Delta^-_i,$ where $\Delta^+_i$ consists of these $(-2)$-twigs of $D_i$ which meet $\Upsilon_i$. It is known that the components of $\Upsilon_i$ are disjoint, so each connected component of $\Delta_i^+$ meets exactly one component of $\Upsilon_i$. We put $$\Bk' \Delta_i^-=\Bk_{D_i}(\Delta_i^-)$$ (see Section \ref{sec:preliminaries}), i.e.\ $\Bk'\Delta_i^-$ is the unique (effective, fractional) $\Q$-divisor supported on $\Supp \Delta_i^-$ such that for every component $R$ of $\Delta_i^-$ we have $(K_i+D_i-\Bk'\Delta_i^-)\cdot R=0$, equivalently $\Bk'\Delta_i^-\cdot R$ equals $-1$ if $R$ is the tip of $D_i$ contained in $\Delta_i^-$ and is $0$ otherwise. One computes \beq \alpha_i^*(K_{Y_i}+\frac{1}{2}D_{Y_i})=K_i+\frac{1}{2}D_i^\flat,\eeq where $$D_i^\flat=D_i-\Upsilon_i-\Delta^+_i-\Bk'\Delta^-_i.$$ Now the log MMP works so that either $(Y_i,\frac{1}{2}D_i)$ is a log Mori fiber space (which in our case is impossible, because of the inequality $\kappa(K+\frac{1}{2}D)\geq 0$) or $K_{Y_i}+\frac{1}{2}D_{Y_i}$ is nef, and then we put $i=n$ (the process stops), or there is an extremal $(K_{Y_i}+\frac{1}{2}D_{Y_i})$-negative contraction $\theta_i\:Y_i\to Z$ which is birational. Consider the latter case. We denote the proper transform of $\Exc \theta_i$ on $X_i$ by $A_i$. We have $(K_i+\frac{1}{2}D_i^\flat)\cdot A_i<0.$ Because of the negativity of the contracted locus the existence of the contraction $\theta_i$ implies (\cite[Corollary 3.5]{Palka-minimal_models}) that \beq\label{eq:Ai} A_i\cdot (\Upsilon_i+\Delta_i^+)=0 \text{\ \ and\ \ } A_i\cdot (D_i-\Delta^-_i)=A_i\cdot\Delta^-_i=1\eeq and that the component of $\Delta^-_i$ meeting $A_i$ is a tip of $\Delta^-_i$ (but not necessarily a tip of $D_i$). What is crucial, $A_i$ is a $(-1)$-curve whose affine part, i.e.\ $A_i\cap (X_i\setminus D_i)$, is isomorphic to $\C^*$. Finally, $\psi_{i+1}\:(X_i,D_i)\to (X_{i+1},D_{i+1}),$ with $D_{i+1}=\psi_{i*}D_i$, is by definition the composition of successive contractions of superfluous $(-1)$-curves in $D_i+A_i$ and its images starting from $A_i$. Because $A_i$ meets two components of $D_i$, the contractions are inner with respect to $D_i+A_i$. It follows that \beq\label{eq:psi_is_inner}\psi_{i+1}^*(K_{i+1}+D_{i+1})=K_i+D_i+A_i.\eeq We say that $\psi_{i+1}$ is of \emph{type $II$} if it contracts both components of $D_i$ meeting $A_i$; otherwise it is of \emph{type $I$}. In any case, all components of $D_i$ contracted by $\psi_{i+1}$ are contained in maximal twigs meeting $A_i$. The morphism $\psi_{i+1}'\:(X_i',D_i')\to (X_{i+1}',D_{i+1}')$ is a lift of $\psi_i$ and the morphism $\psi_{i+1}''\:(Y_i,D_{Y_i})\to (Y_{i+1},D_{Y_{i+1}})$ comes by factorizing $\alpha_{i+1}\circ\psi_{i+1}$ through $\alpha_i$. We put $\psi=\psi_n\circ\ldots\circ\psi_1\:(X_0,D_0)\to (X_n,D_n).$

\bdfn The pairs $(Y_n,\frac{1}{2}D_{Y_n})$ and $(X_n',\frac{1}{2}D_n')$ constructed above are called a \emph{minimal model} and respectively an \emph{almost minimal model}  of $(X_0,\frac{1}{2}D_0)$. \edfn

Let $E$ and $E_i$ for $i\geq 0$ be the proper transforms of $\bar E$ on $X$ and on $X_i$ respectively. By construction $X_i$ and $E_i$ are smooth and $D_i-E_i$ is an snc-divisor. Although $D_i$ is not snc, it has smooth components and contains no superfluous $(-1)$-curves. Recall that by \cite[4.1(viii)]{Palka-minimal_models} the process of minimalization $\psi$ does not affect $h^0(m(2K_i+D_i))$, i.e.\ $$h^0(m(2K_i+D_i))=h^0(m(2K+D)).$$ Also, the behaviour of $\Delta$ and $\Upsilon$ under $\psi$ is well understood. First of all, $\psi_{(i+1)*}(\Upsilon_i)$ is a subdivisor of $\Upsilon_{i+1}$. If $A_i$ meets $\Delta_i^-$ not in a tip of $D_i$ then $\psi_{i+1}$ contracts only $A_i$ (hence is of type $I$) and we have $b_0(\Delta_{i+1})=b_0(\Delta_i)$ and $b_0(\Delta^+_{i+1})=b_0(\Delta^+_i)+1$. If $A_i$ meets $\Delta_i^-$ in a tip of $D_i$ then $\psi_{i+1}$ contracts at least $A_i$ and the connected component of $\Delta_i^-$ meeting $A_i$ (with a minor exception when this connected component and $A_i$ meet the same irreducible component of $D_i-\Delta_i^-$). In the latter case $b_0(\Delta_{i+1})=b_0(\Delta_i)-1$ and either $b_0(\Delta^+_{i+1})=b_0(\Delta^+_i)$ or $b_0(\Delta^+_{i+1})=b_0(\Delta^+_i)+1$.

Denote the cusps of $\bar E$ by $q_1,\ldots, q_c$. We introduce the following numbers characterizing the geometry of the boundary. Assume $j\in \{1,\ldots,c\}$. We write $\tau_j$ for the number of times $\psi_0$ touches $E$. Equivalently, $\tau_j$ is the number of curves over the cusp $q_j$ contracted by $\psi_0$. We put $s_j=1$ if $\psi_0$ contains (in a decomposition into blowdowns) a contraction over $q_j$ which is outer for $D-E$ and $s_j=0$ otherwise. Put $\tau_j^*=\tau_j-s_j-1$ (since $\tau_j\geq 2$ we have $\tau_j^*\geq 0$), $\tau^*=\ds\sum_{j=1}^c \tau_j^*$ and $s=\ds\sum_{j=1}^c s_j$.

\blem\label{lem:(Xi,Di)_properties} Let $(X_i,D_i)$, $(X_i',D_i')$ and $\psi_i$, $i=0,\ldots,n$ be as above. \benum[(i)]

\item $D_i'$ is snc-minimal,

\item $K\cdot (K+D)=p_2(\PP^2,\bar E)$,

\item $\#D_i=\rho(X_i)+i$ for $i\geq 0$,

\item $p_2(\PP^2,\bar E)+i+\ind(D_i')\leq 5$. In particular, $p_2(\PP^2,\bar E)\leq 4$.

\item For every component $U$ of $D_i-\Delta_i$ we have $U\cdot \Delta_i\leq 1$.

\item For every component $U$ of $D_0-E_0$ $$\psi(U)\cdot E_n\leq U\cdot E_0+1.$$ If the equality holds then the unique $\psi_i$ increasing the intersection of the images of $U$ and $E_0$ is of type $I$. Moreover, it touches $U$ exactly once and either $U$ is the component of $\Delta_0^-$ met by $A_{i-1}$ or there is a unique connected component of $\Delta_0^-$ meeting $U$ and this component is contracted by $\psi_i$.
\eenum\elem

\begin{proof} The proof of part (i) given in \cite[3.7]{Palka-minimal_models} for $D_n'$ works for any $i$. For (ii)-(iii), (v) and (vi) see \cite[4.3(i),(ii) and 4.1(vi),(vii)]{Palka-minimal_models} respectively. We note that (ii) is based on the Kawamata-Viehweg vanishing theorem.

(iv) The log resolution $\varphi_n\:(X_i',D_i')\to (X_i,D_i)$ does not touch $A_i$, so $A_i'$, the proper transform of $A_i$ on $X_i'$, is a $(-1)$-curve. The equation \eqref{eq:psi_is_inner} implies $(\psi_{i+1}')^*(K_{i+1}'+D_{i+1}')=K_i'+D_i'+A_i'$, hence $(K_{i+1}'+D_{i+1}')^2=(K_i'+D_i')^2+1$. Then (ii) gives $(K_i'+D_i')^2=(K+D)^2+i=p_2(\PP^2,\bar E)+i-2.$ Put $\cal P_i=(K_i'+D_i')^+$. The surface $X_i\setminus D_i$ contains no $\C^1$, because $X\setminus D$ contains none. By (i) and \ref{lem:Bk}(iii) we obtain $$\cal P_i^2=(K_i'+D_i')^2+\ind (D_i')=p_2(\PP^2,\bar E)+i-2+\ind (D_i').$$ Because $A_i\cap(X_i\setminus D_i)\cong \C^*$, we have $\chi(X_{i+1}\setminus D_{i+1})=\chi(X_i\setminus D_i)$. The logarithmic Bogomolov-Miyaoka-Yau inequality gives $\cal P_i^2\leq 3\chi(X_i'\setminus D_i')=3\chi(X\setminus D)=3$.
\end{proof}

We say the cusp $q_j\in \bar E$ is \emph{semi-ordinary} if it is locally analytically isomorphic to the singular point of $x^2=y^{2m+1}$ at $0\in\Spec \C[x,y]$ for some $m\geq 1$. Let $c_0$ and $c_1$ be the numbers of semi-ordinary and non-semi-ordinary cusps of $\bar E$ respectively. Note that $\#\Upsilon_0=c_0$. Put $\eta=\#\Upsilon_n-\#\Upsilon_0$. For $k=0,1$ let $n_k$ be the number of contracted $A_i$'s, i.e.\ the $(-1)$-curves defined above, for which $A_i\cdot E_i=k$. By \eqref{eq:Ai} $n=n_0+n_1$. By \cite[4.3]{Palka-minimal_models}) \beq\label{eq:Kn(Kn+Dn)} K_n\cdot (K_n+D_n)=p_2(\PP^2,\bar E)-c-\tau^*-n \eeq and \beq\label{eq:En(Kn+Dn)} E_n\cdot (K_n+D_n)=2c-2+\tau^*+n_1.\eeq

\subsection{The log MMP and the Cremona equivalence} The assumption that $\bar E\subseteq \PP^2$ is not Cremona equivalent to a line, or equivalently that $2K_X+E\geq 0$, was not used in the construction above. Now it will allow us to analyze $D_n$ and the process of minimalization in more detail. For example if we intersect $2K_n+E_n$, which is effective, with $\alpha_n^*(K_{Y_n}+\frac{1}{2}D_{Y_n})$, which is nef, we get a non-negative number (see \ref{prop:(Xi,Di)_properties_for_CN}(i) for the resulting inequality).

Let $C_j$ be the $(-1)$-curve of $D_0-E_0$ over $q_j$ and let $\cal C=C_1+\ldots+C_c$. All $C_j$ are tangent to $E_0$, so are not contracted by $\psi$. Denote their images on $X_n$ by $C_1',\ldots,C_c'$. Let $\cal C_+$ and $\cal C_{exc}$ be the sums of these $C_j'$'s whose self-intersection is non-negative and equal to $(-1)$ respectively. Clearly, $\psi_*(\cal C)=\cal C_++\cal C_{exc}$. In the process of minimalization $\psi$ some new non-superfluous $(-1)$-curves might have been created in the boundary. Let $\cal L$ be the sum of $(-1)$-curves of $D_n-E_n$ which are not components of $\cal C_{exc}$. Put $n_{exc}=\#\cal L$. We say that a $(-1)$-curve $L$ in $D_i-E_i$ (and its image in $D_n-E_n$, if it belongs to $\cal L$) is \emph{created by $\psi_i$} if and only if $\psi_{i*}^{-1}(L)$ is not a $(-1)$-curve. Put $$R_n=D_n-E_n-\cal C_+-\cal C_{exc}-\cal L.$$

Note that if $L$ is a common component of $\cal L$ and $\Upsilon_n$ then $L'=\psi_*^{-1}(L)$ is not a component of $\cal C$, which implies that $L\cdot E_n\leq 1$. Indeed, otherwise \ref{lem:(Xi,Di)_properties} gives $L'\cdot E_0\neq 0$, so $L'\cdot E_0=L'\cdot C_j=1$ for some $j$ and hence $L\cdot (D_n-\Delta_n)\geq L\cdot (C_j'+E_n)\geq 3$, which contradicts the definition of $\Upsilon_n$. Therefore, we may decompose $\eta=\#\Upsilon_n-\#\Upsilon_0$, introduced above, as $\eta=\eta_0+\eta_1$, where for $k=0,1$ the number $\eta_{k}$ counts the common components of $\cal L$ and $\Upsilon_n$ which meet $E_n$ exactly $k$ times.

\blem\label{lem:properties_of_(Xi,Di)_for_CN} Assume $1\leq i\leq n$ and $1\leq j\leq c$. The following hold. \benum[(i)]

\item Let $L\subseteq X_n$ be the image of a $(-1)$-curve in $D_i-E_i$ created by some $\psi_i$, $i\geq 1$ and let $M$ be a component of $D_n-E_n-L$ with $m=-M^2$ and $u=L\cdot M>0$. Then $L\leq \cal L$ and if $m\leq u^2$ then $M\cdot E_n+2m\geq 4+(2-L\cdot E_n) u.$

\item Each $\psi_i$ creates at most one $(-1)$-curve.

\item The components of $\cal L$ are disjoint and $n_{exc}\leq n\leq 4$.

\item For $k=0,1$ we have $\eta_k\leq n_k$.

\item Components of $R_n$ intersect non-negatively with $K_n$.

\item If $(C_j')^2\geq 0$ then $K_n\cdot C_j'+\tau^*_j\geq 0$.

\item Assume $\psi_i$ creates a $(-1)$-curve $L$ in $D_i-E_i$. Then for every component $V$ of $D_i-E_i-L$ we have $V\cdot E_i=\psi_{i*}^{-1}(V)\cdot E_{i-1}$. If $L\cdot E_i\neq \psi_{i*}^{-1}(L)\cdot E_0$ then $A_{i-1}\cdot E_{i-1}=1$, $L\cdot E_n=\psi_{*}^{-1}(L)\cdot E_0+1$ and $\psi_i$ contracts either $A_{i-1}$ or $A_{i-1}$ together with a $(-2)$-twig of $D_{i-1}$ meeting $\psi_{i*}^{-1}(L)$.
\eenum \elem

\begin{proof}
(i) Note that once $\psi_i$ creates a $(-1)$-curve $L$, this curve is not touched by $\psi_j$ for $j>i$. Indeed, otherwise by \ref{lem:(Xi,Di)_properties}(vi) the first blowdown touching $L$ makes it into a $0$-curve intersecting the image of $E_0$ at most twice, which contradicts \ref{lem:cuspidal_of_gt_and_khalf}(iv). Put $\mu =2-L\cdot E_n$ and $f=uL+M$. By \ref{lem:(Xi,Di)_properties}(vi) $\mu \geq 0$. We have $f\cdot L=0$ and $f\cdot M=u^2-m\geq 0$, so $f$ is nef.  Because $2K_0+E_0\geq 0$, we have $2K_n+E_n\geq0$. But $\mu L$ is in the fixed part of $|2K_n+E_n|$, so $0\leq f\cdot (2K_n+E_n-\mu L)=M\cdot (2K_n+E_n-\mu L)=2(m-2)+M\cdot E_n-\mu u$, which is equivalent to the inequality (i).

(ii)-(iii) Since for $i\geq 1$ the morphism $\psi_i$ is a composition of blowdowns which are inner for $D_{i-1}+A_{i-1}$, we see that it can create at most two $(-1)$-curves. Suppose two components $L_1, L_2$ of $\cal L$ meet. (This is for example the case if some $\psi_i$ creates more than one $(-1)$-curve.) For $j=1,2$ let $L_j'$ be the proper transform of $L_j$ on $X_0$. Since $L_j'$ is not a $(-1)$-curve, we have $L_j'\cdot E_0\leq 1$, so by \ref{lem:(Xi,Di)_properties}(vi) $L_j\cdot E_n\leq L_j'\cdot E_0+1\leq 2$. Since $L_1+L_2$ is nef, by \ref{lem:cuspidal_of_gt_and_khalf}(iii) $4\leq (L_1+L_2)\cdot E_n$, which implies $L_1\cdot E_n=L_2\cdot E_n=2$ and hence $L_1'\cdot E_0=L_2'\cdot E_0=1$. Therefore $L_j'$'s lie over different cusps. It follows also that for $j=1,2$ there exists $\psi_{m_j}$, which increases the intersection of the image of $L_j'$ with the image of $E_0$. But then $A_{m_1-1}$ and $A_{m_2-1}$, and hence their proper transforms on $\pi(X)=\PP^2$ are disjoint; a contradiction. Thus each $\psi_i$ creates at most one $(-1)$-curve and these $(-1)$-curves remain disjoint on $X_n$. We have $n\leq 5-p_2(\PP^2,\bar E)\leq 4$ by \ref{lem:(Xi,Di)_properties}(iv).

(iv) Fix $i\geq 1$. Assume $U$ is a component of $D_{i-1}-E_{i-1}-\Upsilon_{i-1}$ such that $\psi_i(U)$ is a component of $\Upsilon_i$. By (ii) it is enough to show that $A_{i-1}$ meets $E_{i-1}$ if and only if $\psi_i(U)$ meets $E_i$. Since $\psi_{j*}\Upsilon_{j-1}\leq \Upsilon_j$, the image of $U$ on $X_n$ is a component of $\Upsilon_n$ and in particular $\psi_i(U)$ is not touched by $\psi_j$ with $j>i$. If $A_{i-1}\cdot E_{i-1}=1$ then $\psi_i$ contracts either only $A_i$ or $A_i$ and a connected component of $\Delta_{i-1}^-$ meeting it. In both cases $\psi_i(U)$ meets $E_i$. Suppose $A_{i-1}\cdot E_{i-1}=0$ and $\psi_i(U)\cdot E_i>0$. Since $\psi_i(U)$ is a component of $\Upsilon_i$, the proper transform $U_0$ of $U$ on $X_0$ does not meet $E_0$. Indeed, otherwise it would also meet some $C_j$ and a $(-2)$-twig of $D_0$ not touched by $\psi$, which implies $\beta_{D_i}(U)>3$, in contradiction to the definition of $\Upsilon_i$. Thus some $\psi_j$ with $j<i$ increases the intersection of images of $U_0$ and $E_0$, so by \ref{lem:(Xi,Di)_properties}(vi) $\psi_j$ contracts some $(-2)$-twig of $D_{j-1}$ meeting $U_0$. Note that the image of $U_0$ on $X_j$ is not a $(-1)$-curve. Because $\psi_i(U)$ meets $E_i$ transversally and is not superfluous, we have $\beta_{D_i}(\psi_i(U))\geq 3$. Because $\psi_i(U)$ is a component of $\Upsilon_i$, it follows that $U_0$ meets more than one $(-2)$-twig of $D_0$; a contradiction with \ref{lem:(Xi,Di)_properties}(v).

(v)-(vi) Let $L$ be a component of $R_n+\cal C_+$ with $K_n\cdot L<0$. By definition $L$ is not a $(-1)$-curve, so $L^2\geq 0$. We blow up $L^2$ times on $L$ and we denote the proper transforms of $L$ and $E$ by $L'$ and $E'$ respectively. Now $L'$ induces a $\PP^1$-fibration of the blowup, and by \ref{lem:cuspidal_of_gt_and_khalf}(iii) $L'\cdot E'\geq 4$. By \ref{lem:(Xi,Di)_properties} this can happen only if $\psi_*^{-1}(L)\cdot E_0\geq 3$, so $L=C_j'$ for some $j$, say $j=1$. We get $4\leq L'\cdot E'=L\cdot E_n-L^2\leq C_1\cdot E_0+1-L^2,$ so $\tau_1^*+K_n\cdot C_1'\geq -s_1$. Suppose $\tau_1^*+K_n\cdot C_1'<0$. Then $s_1=1$ and all the inequalities become equalities, so $\tau_1^*+K_n\cdot C_1'=-1$ and $C_1'\cdot E_n=C_1\cdot E_0+1$, which implies that $C_1$ meets some $(-2)$-twig of $D_0$. Since $s_1=1$, $C_1$ together with this $(-2)$-twig is a connected component of $D_0-E_0$. Moreover, $\psi$ touches $C_1$ once, exactly when contracting this $(-2)$-twig together with some $A_{i_0}$ intersecting it in a tip. But this $A_{i_0}$ does not meet the exceptional divisor over cusps other than $q_1$, so $\pi_0(A_{i_0})$ is a $0$-curve on $\PP^2$; a contradiction.

(vii) Assume $V\cdot E_i>\psi_{i*}^{-1}(V)\cdot E_{i-1}$ for some component $V$ of $D_i-E_i$. By \ref{lem:(Xi,Di)_properties}(ii) $A_{i-1}\cdot E_{i-1}=1$, $\psi_i$ is of type $I$ and $V=L$, so $\psi_i$ touches $L'=\psi_{i*}^{-1}(L)$ once. Furthermore, either $L'$ is a part of a $(-2)$-twig of $D_{i-1}$ and $A_{i-1}$ meets $L'$ or $L'$ and $A_{i-1}$ meet the same maximal $(-2)$-twig of $D_{i-1}$. By \ref{lem:(Xi,Di)_properties}(v) for a given $V$ this increase of the intersection with the image of $E_0$ can happen for at most one $\psi_i$.
\end{proof}

The proof of Theorem \ref{thm:MAIN1}, given in the next section, is rather lengthy and requires ruling out some very concrete types of cuspidal curves. The basic role in reducing the proof to these cases is played by the following proposition. For $i=0,\ldots,n$ put $\gamma_i=-E_i^2$ and $\gamma=-E^2$. Recall that by \ref{lem:cuspidal_of_gt_and_khalf}(iii) $2K_i+E_i\geq 0$, so $\gamma_i\geq 4$.

\bprop\label{prop:(Xi,Di)_properties_for_CN} Let the notation be as above. Put $\zeta=K_n\cdot (K_n+E_n)$. Then:\ \benum[(i)]

\item $\gamma_n+\tau^*+(n_1-\eta_1)+2(n_0-\eta_0)\leq 2\zeta+2p_2(\PP^2,\bar E)$,

\item $$\sum_{j:C_j'\leq \cal C_+}(K_n\cdot C_j'+\tau_j^*)+\#\cal C_++\sum_{j:C_j\leq \cal C_{exc}}\tau_j^*+(n-n_{exc})+K_n\cdot R_n=p_2(\PP^2,\bar E)-\zeta.$$
\eenum\eprop

\begin{proof} By construction the divisor $2K_n+D_n^\flat$ is nef, so $$0\leq (2K_n+D_n^\flat)\cdot(2K_n+E_n)=(2K_n+D_n)\cdot(2K_n+E_n)-\Upsilon_n\cdot(2K_n+E_n).$$ We have $$\Upsilon_n\cdot (2K_n+E_n) =\Upsilon_0\cdot(2K_0+E_0)-2\eta_0-\eta_1=-2\eta_0-\eta_1$$ and $$(2K_n+D_n)\cdot(2K_n+E_n)=2\zeta-K_n\cdot E_n+2K_n\cdot(K_n+D_n)+E_n\cdot (K_n+D_n).$$ By \eqref{eq:Kn(Kn+Dn)} and \eqref{eq:En(Kn+Dn)} the latter expression equals $2\zeta-\gamma_n+2p_2(\PP^2,\bar E)-\tau^*-n_1-2n_0$, which proves (i). We have $$K_n\cdot \sum_{j}C_j'=K_n\cdot \cal C_+-(c-\#\cal C_+)\text{\ \ and\ \ } K_n\cdot \cal L=-n_{exc}.$$ By \eqref{eq:Kn(Kn+Dn)} $K_n\cdot (D_n-E_n)+\zeta+\tau^*+c+n=p_2(\PP^2,\bar E)$, which leads to (ii).
\end{proof}

Let us now transfer some results to the level of $(X_n',D_n')$. Let $\cal U$ be the sum of the $(-1)$-curves of $D$. Note that these are exactly the components of $D-E$ meeting $E$. By the definition of $\psi$ and $\varphi_n$, they are not touched by $\psi'\:(X,D)\to (X_n',D_n')$, i.e.\ their images are $(-1)$-curves in $D_n'$. Let $\cal L'$ be the sum of $(-1)$-curves in $D_n'-E_n'-\psi'_*(\cal U)$, put $n_{exc}'=\#\cal L'$. As in \ref{lem:properties_of_(Xi,Di)_for_CN}(iii) we have $n_{exc}'\leq n$. Put $R_n'=D_n'-E_n'-\psi'_*(\cal U)-\cal L'$.

\bcor\label{cor:b2_bounds_etc} Let the notation be as above. The almost minimal model $(X_n',\frac{1}{2}D_n')$ of $(X_0,\frac{1}{2}D_0)$ has the following properties: \benum[(i)]

\item Components of $R_n'$ intersect $K_n'$ non-negatively,

\item $K_n'\cdot (K_n'+E_n')=K_n\cdot(K_n+E_n)=\zeta$,

\item $K_n'\cdot R_n'=p_2(\PP^2,\bar E)-\zeta-n+n_{exc}'+c$,

\item $2\leq \frac{1}{2}(\gamma_n+\tau^*)\leq\zeta+p_2(\PP^2,\bar E)\leq 2p_2(\PP^2,\bar E)$,

\item $\rho(X_n')\leq 8+2p_2(\PP^2,\bar E)+c+s+\zeta \leq 20+2c$,

\item $\#D_n'\leq 8+2p_2(\PP^2,\bar E)+c+s+\zeta+n\leq 24+2c$.
\eenum\ecor

\begin{proof} (i) By definition $R_n'$ contains no $(-1)$-curves. Let $W$ be a component of $R_n'$ with $W^2\geq 0$. Lemma \ref{lem:(Xi,Di)_properties}(vi) implies that $W\cdot E_n'\leq (\psi')_*^{-1}(W)\cdot E+1=1$, so after blowing up on $W$ until $W^2=0$ we get a $\PP^1$-ruled surface for which the proper transform of $E$ meets a general fiber at most once. This contradicts \ref{lem:cuspidal_of_gt_and_khalf}(iv).

(ii) The centers of all blowups constituting $\varphi_n$ belong to respective proper transforms of $E$, so $\varphi_n^*(K_n+E_n)=K_n'+E_n'$, hence $K_n'\cdot (K_n'+E_n')=K_n\cdot(K_n+E_n)=\zeta.$

(iii) We have $K_n'\cdot (K_n'+D_n')=K\cdot (K+D)-n=p_2(\PP^2,\bar E)-n$ by \ref{lem:(Xi,Di)_properties}(ii), hence by (ii) $K_n'\cdot R_n'=p_2(\PP^2,\bar E)-n+n_{exc}'+c-\zeta$.

(iv)  By \ref{prop:(Xi,Di)_properties_for_CN}(i) $\zeta+p_2(\PP^2,\bar E)\geq \frac{1}{2}(\gamma_n+\tau^*)$. By (ii) and by the Riemann-Roch theorem $\zeta\leq h^0(2K_n'+E_n')=h^0(2K+E)\leq p_2(\PP^2,\bar E).$

(v) We have $K_n'\cdot E_n'=K_n\cdot E_n+\tau=\gamma_n-2+\tau$, so by the Noether formula $\zeta=(K_n')^2+K_n'\cdot E_n'=8-\rho(X_n')+\gamma_n+\tau$. Now \ref{prop:(Xi,Di)_properties_for_CN}(i) gives $\zeta+8+2p_2(\PP^2,\bar E)-\rho(X_n')\geq -c-s$.

(vi) follows from (v) and from the equality $\#D_n'=\rho(X_n')+n$ (see \ref{lem:(Xi,Di)_properties}(iii)).
\end{proof}

For $j=1,\ldots,c$ let $Q_j$ be the reduced exceptional divisor of $\pi\:X\to \PP^2$ over the cusp $q_j\in \bar E$ and let $\wt Q_j=\psi_0(Q_j)$.

\blem\label{lem:E_with_mult<=3_is_rectifiable} Assume the cusps $q_2,\ldots, q_c\in \bar E$ have multiplicity two (equivalently, they are semi-ordinary). Then $q_1$ has multiplicity at least four. \elem

\begin{proof} Let $\mu(q_j)$ denote the multiplicity of $q_j\in \bar E$. If $j\geq 2$ then, because $\mu(q_j)=2$, the divisor $\wt Q_j$ is a chain of type $[(2)_{t_j-1},1]$ for some $t_j\geq 1$. Suppose $\mu(q_1)\leq 3$.

\bcl $\mu(q_1)\neq 2$. \ecl

Suppose $\mu(q_1)=2$. If $n\neq 0$ then $D_0-E_0+A_0$ contains a chain of type $f=[1,2,\ldots,2,1]$ meeting $E_0$ at most three times. But $f$ is nef, so by \ref{lem:cuspidal_of_gt_and_khalf}(iii) $0\leq f\cdot (2K_0+E_0)=-4+f\cdot E_0$; a contradiction. Thus $n=0$ and so \ref{prop:(Xi,Di)_properties_for_CN} gives $K_0\cdot(K_0+E_0)=p_2(\PP^2,\bar E)$ and $\gamma_0\leq 4p_2(\PP^2,\bar E)$. By \ref{lem:(Xi,Di)_properties}(iv) $p_2(\PP^2,\bar E)\leq 4$. The Noether formula gives $\#D_0=\rho(X_0)=10-K_0^2=8+\gamma_0-p_2(\PP^2,\bar E).$ For the resolution $X_0\to \PP^2$ all multiplicities of singular points of $\bar E$ are equal to $2$, so the genus-degree formula gives $$\frac{1}{2}(\deg \bar E-1)(\deg \bar E-2)=\#(D_0-E_0)=7+\gamma_0-p_2(\PP^2,\bar E)\leq 19.$$ The solutions are $(\deg \bar E,\#(D_0-E_0))=(6,10)$ and $(7,15)$. These configurations have been ruled out by Yoshihara \cite{Yoshihara} by analyzing cohomology of double covers of $\PP^2$ branched along $\bar E$ and $\bar E\cup (\text{tangent to } \bar E)$ respectively. In fact we could refer here to the inequality of Matsuoka-Sakai \ref{lem:cuspidal_of_gt_and_khalf}(vi) saying that $\deg\bar E<3\max_j\mu(q_j)=6$, but in our two cases (called A-I and A-II in \cite{MaSa-cusp}) the proof goes by referring to Yoshihara.

\med Therefore $\mu(q_1)=3$. It follows that $\wt Q_1$ is a chain of type $[(2)_{t_1+1-s_1},1]$ for some $t_1\geq 0$ and that $\tau_1=2+s_1$. If $s_1=1$ then $C_1$ is the only component of $\wt Q_1$ meeting $E_0$. If $s_1=0$ then we denote the component of $\wt Q_1-C_1$ meeting $E_0$ by $U$, otherwise we put $U=0$. We have $U\cdot E_0=U\cdot C_1=1$.

\bcl $n=0$. \ecl

Suppose $n\neq 0$. Because $\tau^*=1$, by \ref{lem:properties_of_(Xi,Di)_for_CN}(vi) no $C_j$ is touched by $\psi$. If $A_0$ does not meet $E_0$ or if it meets $D_0-E_0-\wt Q_1$ then $D_0-E_0+A_0$ contains a nef chain $f=[1,2,\ldots,2,1]$ containing $A_0$ and some $C_j$ ($j\geq 2$ in the second case) such that $f\cdot (2K_0+E_0)={f\cdot E_0}-4<0$, which contradicts the fact that $2K_0+E_0\geq 0$. It follows that $A_0$ meets $E_0$ and some tip of the connected component of $\Delta_0$ contained in $\wt Q_1$. If $s_1=1$ then the contraction of $D_0-E_0-C_1+A_0$ maps $X_0$ onto $\PP^2$ and $E_0$ to a cuspidal curve with cusps of multiplicity two. If $s_1=0$ then the same happens for the contraction of $D_0-E_0-U+A_0$ (note $A_0$ does not meet $U$, because $U$ is not a component of $\Delta_0$). But the latter cuspidal curve is Cremona equivalent to $\bar E$, hence Cremona non-equivalent to a line, which is impossible by the previous claim; a contradiction. By \eqref{eq:Kn(Kn+Dn)} $\zeta=p_2(\PP^2,\bar E)-1$.

\bcl $p_2(\PP^2,\bar E)=c=2$ and $\gamma_0=5$. \ecl

By \ref{rem:deg<=5_rectifiable} and \ref{lem:cuspidal_of_gt_and_khalf}(vi) $\deg \bar E\in \{6,7,8\}$. The sequence of characteristic pairs for $Q_1$ (see \ref{lem:ind_using_pairs}) is $\binom{3t_1+5-s_1}{3}$, $t_1\geq 0$ and for $Q_j$, $2\leq j\leq c$ it is $\binom{2t_j+1}{2}$, where $t_j\geq 1$. The Noether formula for $X_0$ gives $\#D_0=\rho(X_0)=10-K_0^2=9+\gamma_0-p_2(\PP^2,\bar E)$, so $t_1+t_2+\ldots+t_c=6+\gamma_0-p_2(\PP^2,\bar E)+s_1$. We have also $\gamma=\gamma_0+2c+s_1$. The equations in \ref{cor:I and II} give:
\begin{align*}
3d&=t_1+19+\gamma_0-2p_2(\PP^2,\bar E),\\
d^2&=5t_1+37+3\gamma_0-4p_2(\PP^2,\bar E).
\end{align*}

By \ref{prop:(Xi,Di)_properties_for_CN}(i) $4\leq\gamma_0\leq 4p_2(\PP^2,\bar E)-3$. The solutions are:
$(\gamma_0,d,p_2(\PP^2,\bar E),t_1)=(5,7,2,1)$, $(5,8,2,4)$, $(8,7,3,0)$, $(8,8,3,3)$, $(11,8,4,2)$. Using \ref{lem:ind_using_pairs} we compute $$\ind (D)=1-\frac{3}{3t_1+5-s_1}+\frac{1+s_1}{3}+\sum_{j=2}^c(\frac{3}{2}-\frac{2}{2t_j+1})\geq \frac{11}{15}+(c-1)\frac{5}{6}.$$ The Bogomolov-Miyaoka-Yau inequality \ref{lem:(Xi,Di)_properties}(iv) gives $\ind (D)\leq 5-p_2(\PP^2,\bar E)\leq 3$. It follows that $c\leq 3$. Only the first two solutions satisfy the BMY inequality and only if $c=2$. We obtain $\gamma_0=5$, $p_2(\PP^2,\bar E)=2$, $t_1+t_2=9+s_1$ and $(d,t_1)=(7,1)$ or $(8,4)$.

\bcl There is a $(-1)$-curve $V$ on $X_0$, such that $2K_0+E_0\sim V$. \ecl

Write $2K_0+E_0\sim V_1+\ldots +V_k$, where $V_i$, $i\leq k$ are irreducible. We have $n=0$, $\tau^*=1$ and $K_0\cdot (D_0-E_0-\cal C)=0$, so by \ref{prop:(Xi,Di)_properties_for_CN}(ii) $\zeta=p_2(\PP^2,\bar E)-1=1$. Hence $(2K_0+E_0)^2=4\zeta+E_0^2=-1$, so there is a component $V$ among $V_i$'s, such that $V\cdot (2K_0+E_0)<0$. It follows that $V^2<0$ and $V\cdot K_0<0$, so $V$ is a $(-1)$-curve for which $V\cdot E_0\leq 1$. In particular, $V\not\subseteq D_0$.

We now show that $V\cdot (\Delta_0^++\Upsilon_0)=0$, $V\cdot \Delta_0^-\leq 1$ and $V\cdot U\leq 1$. Note that $\Upsilon_0=C_2$. Suppose $V$ meets $\Delta_0^++C_2$. Let $f$ be the shortest subchain of $\Delta_0^++C_2$ containing $C_2$ and a component of $\Delta_0^++C_2$ meeting $V$. Then $f+V$ is nef and $$0\leq (f+V)\cdot (2K_0+E_0)=V\cdot E_0-2<0;$$ a contradiction. Suppose $V\cdot \Delta_0^-\geq 2$. Let $f=U_1+U_2+\ldots+U_m$ be a subchain of $\Delta_0^-$, such that $U_1\cdot V$, $U_m\cdot V>0$ and $U_1\cdot V\geq 2$ if $m=1$. Then $f+V$ is nef and $(f+V)\cdot (2K_0+E_0)\leq -2+V\cdot E_0<0$; a contradiction. Suppose $V\cdot U\geq 2$. Then $f=U+2V$ is nef and $0\leq f\cdot (2K_0+E_0)\leq 2V\cdot E_0-3$; a contradiction.

Because $\gamma_0=5$, we get $K_0^2=-2$ and then (see the proof of \ref{prop:(Xi,Di)_properties_for_CN}(i)) $$(2K_0+D_0^\flat)\cdot(2K_0+E_0)=2\zeta+2p_2(\PP^2,\bar E)-\gamma_0-\tau^*=0.$$ Since $(2K_0+D_0^\flat)$ is nef, we get $(2K_0+D_0^\flat)\cdot V=0$, hence $(2K_0+D_0)\cdot V=\Bk' \Delta_0^-\cdot V$. But $\Delta_0^-\cdot V\leq 1$ and $\Bk' \Delta_0^-$ has proper fractional coefficients, so $\Delta_0^-\cdot V=0$. It follows that $(D-E_0-C_1-U)\cdot V=0$ and $(E_0+C_1+U)\cdot V=2$. If $V$ does not meet $C_1$ then $V\cdot E_0=V\cdot U=1$. If $V$ meets $C_1$ then $f=V+C_1$ is nef, hence $0\leq f\cdot (2K_0+E_0)=-4+\tau_1+V\cdot E_0$, so $V\cdot E_0=1$ and $\tau_1=3$. In total, the above computations imply that the effective divisor $2K_0+E_0-V$ intersects trivially with all components of $D_0$. Because $\Pic X_0\otimes \Q$ is generated by the components of $D_0$, we obtain $2K_0+E_0-V\sim 0$. Note also that $V\cap (X_0\setminus D_0)$ is isomorphic to $\C^*$.

\med The proper transform $V'$ of $V$ on $X$ is a $(-1)$-curve. The surface $X\setminus(D\cup V')$ contains no $\C^1$'s, because $X\setminus D$ contains none. Let $(X,D)\to (Y,D_Y)$ be the contraction of $V'$. We have $(K_Y+D_Y)^2=(K+D)^2+1=p_2(\PP^2,\bar E)-1$, so since $(Y,D_Y)$ is almost minimal, $$((K_Y+D_Y)^+)^2=p_2(\PP^2,\bar E)+\ind (D_Y)-1=\ind (D_Y)+1.$$ We compute $$\ind (D_Y)=1-\frac{1}{t_1+1}+\frac{1+s_1}{3}+(\frac{3}{2}-\frac{2}{2t_2+1}).$$ The Bogomolov-Miyaoka-Yau inequality for $(Y,D_Y)$ reads as $$((K_Y+D_Y)^+)^2\leq 3\chi(X\setminus(D\cup V'))=3.$$ Thus $\ind (D_Y)\leq 2$. In both cases $(t_1,t_2)=(1,8+s_1)$ and $(4,5+s_1)$ this inequality fails; a contradiction.
\end{proof}

\section{Proof of Theorem \ref{thm:MAIN1}}

With the description of the minimalization process from Section \ref{sec:logMMP} we are now ready to prove Theorem \ref{thm:MAIN1}. The crucial role is played by the inequality \ref{prop:(Xi,Di)_properties_for_CN}(i) in which small $p_2(\PP^2,\bar E)$ gives bigger restrictions on the minimal model and hence on the cuspidal curve $\bar E\subseteq \PP^2$, eventually leading to a contradiction for $p_2(\PP^2,\bar E)\leq 2$. We keep the assumptions and notation from the previous section. Recall that $D_n=R_n+E_n+\cal C_++\cal C_{exc}+\cal L$. Put $R_0=D_0-E_0-\cal C$.

\blem\label{prop:jump>=2}  We have $p_2(\PP^2,\bar E)\geq K_n\cdot (K_n+E_n)+1$. If the equality holds then $n=1$, $n_{exc}=n_1=0$, $p_2(\PP^2,\bar E)=2$ and $2K_1+E_1\sim 0$. \elem

\begin{proof} Suppose $p_2(\PP^2,\bar E)\leq \zeta+1$. By \ref{lem:properties_of_(Xi,Di)_for_CN}(iii)-(iv) all summands on the left hand side of \ref{prop:(Xi,Di)_properties_for_CN}(ii) are non-negative, so $\zeta=K_n\cdot(K_n+E_n)\leq p_2(\PP^2,\bar E)$ (in fact the Riemann-Roch theorem gives $\zeta\leq h^0(2K_n+E_n)$). Since $\gamma_n\geq 4$, \ref{prop:(Xi,Di)_properties_for_CN}(i) gives $\zeta\geq 1$.

\bcl The lemma holds in case $n\neq n_{exc}$. \ecl

By \ref{lem:properties_of_(Xi,Di)_for_CN}(iii) $n_{exc}\leq n$. Suppose $n\neq n_{exc}$. By \ref{prop:(Xi,Di)_properties_for_CN}(ii) $p_2(\PP^2,\bar E)=\zeta+1$, $n-n_{exc}=1$ and $\#\cal C_+=\tau^*=K_n\cdot R_n=0$. In particular, $R_n$ consists of $(-2)$-curves. Suppose $n_{exc}\neq 0$. There exist a component $L$ of $\cal L$, a component $L'$ of $\cal C_{exc}+\cal L-L$ and components $U_1,\ldots, U_k$ of $R_n$, such that $f=L'+U_1+\ldots+U_k+L$ is connected and has no $(-2)$-tips, hence is nef. Because $\tau_1^*=0$, \ref{lem:properties_of_(Xi,Di)_for_CN}(vii) gives $L\cdot E_n\leq \psi_*^{-1}(L)\cdot E_0+1=1$. We have $0\leq f\cdot (2K_n+E_n)=-4+f\cdot E_n$, so $f\cdot E_n\geq 4$. In particular $E_n\cdot (U_1+\ldots +U_k)\geq 4-E_n\cdot (L+L')\geq 1$, because $L'\cdot E_n\leq 2$ if $L'\leq \cal C_{exc}$ and $L'\cdot E_n\leq 1$ otherwise. Thus among $U_i$'s there is a component $U$ meeting $E_n$. Since $\tau^*=0$, $U_0=\psi_*^{-1}(U)$ is disjoint from $E_0$. Let ${i_0}$ be the smallest $i>0$, such that the images of $U_0$ and $E_0$ on $X_{i}$ meet. By \ref{lem:properties_of_(Xi,Di)_for_CN}(vii) $U\cdot E_n=1$ and $\psi_{i_0}$ creates no $(-1)$-curves in $D_i-E_i$. Because $n-n_{exc}=1$, $U$ is the only component of $R_n$ meeting $E_n$ and every other $\psi_i$ creates a $(-1)$-curve. Say $U_0\leq \wt Q_1$. We obtain $f\cdot (L+L')\geq 3$, hence $L\cdot E_n=1$ and $L'=C_1$. Because $U$ is contained in every chain $f$ as above we infer that $\psi_*^{-1}(\cal L)\leq \wt Q_1$ and hence that in fact $\cal L=L$. It follows that $n=2$ and that the cusps $q_2,\ldots, q_c$ are semi-ordinary. We have $\psi_*^{-1}(L)\cdot E_0<L\cdot E_n$ and $U_0\cdot E_0< U\cdot E_n$, so \ref{lem:(Xi,Di)_properties}(vi) implies that $U_0^2=-3$ and $R_0-U_0$ consists of $(-2)$-curves. Then $\wt Q_1$ is a fork with a branching $(-2)$-curve and maximal twigs of types $[2]$, $[(2)_{t_1},3]$ and $[1,(2)_{t_2}]$ for some $t_1\geq 1$ and $t_2\geq 0$ (cf. \ref{lem:Q_with_small_KQ}(ii); by convention the last curve in the chain meets the branching component) and $\psi_{i_0}$ contracts $[(2)_{t_1}]$. But then $U$ is a tip of $D_n-E_n$, so it is not a part of $U_1+\ldots+U_k$ as above; a contradiction.

Therefore, we have $n_{exc}=0$ and $n=1$. Suppose that $A_0$ meets $E_0$. Then $\psi$ contracts exactly $A_0$ together with some $(-2)$-twig of $D_0$, so $R_0$ consists of $(-2)$-curves and one $(-3)$-curve $V$. This again implies that $\wt Q_1$ is a fork with a branching $(-2)$-curve and maximal twigs of types $[2]$, $[(2)_{t_1},3]$, $[1,(2)_{t_2}]$ for some $t_1,t_2 \geq 0$. Furthermore, $A_0$ touches the tip of $D_0$ contained in $[(2)_{t_1}]$. The contraction of $A_0+D_0-E_0-V$ maps $X_0$ to onto $\PP^2$ and $E_0$ onto a cuspidal curve with only semi-ordinary cusps. By \ref{lem:E_with_mult<=3_is_rectifiable} this is a contradiction. Thus $A_0\cdot E_0=0$. It follows that components of $R_1$ intersect $2K_1+E_1$ trivially. Because $\tau^*=\#\cal C_+=0$, in fact all components of $D_1-E_1$ intersect $2K_1+E_1$ trivially. But $X\setminus D$ is $\Q$-acyclic, so the components of $D$ generate $\NS(X)\otimes \Q$, hence the components of $D_1$ generate $\NS(X_1)\otimes \Q$. We have $2K_1+E_1\geq 0$, so $2K_1+E_1\sim 0$ if and only if $E_1\cdot (2K_1+E_1)=0$. Now in case $p_2(\PP^2,\bar E)=2$ we have $\zeta=1$, so \ref{prop:(Xi,Di)_properties_for_CN}(i) gives $\gamma_1=4$ and hence $E_1\cdot (2K_1+E_1)=E_1\cdot K_1-2=0$, so we are done.

We may therefore assume that $p_2(\PP^2,\bar E)\geq 3$. By \ref{lem:(Xi,Di)_properties}(iv) $\ind(D_1')\leq 1$. Since $\tau^*=0$, $\psi_0$ contracts exactly $c$ $(-2)$-tips of $D_1$, so $\ind(D_1')\geq \frac{1}{2}c$, hence $c\leq 2$. If $c=2$ then the contribution from other twigs to $\ind(D_1')$ is zero, so both $\wt Q_j$, $j=1,2$ are chains and hence are of type $[(2)_{t_j},1]$ for some $t_j\geq 0$. But in the latter case $\wt Q_j-C_j$ is a part of $\Delta_0^+$, hence is not touched by $\psi$; a contradiction. Thus $c=1$. Then the contribution to $\ind(D_1')$ from $E_1$ and the $(-2)$-tip contracted by $\psi_0$ is $\frac{1}{2}+\frac{1}{\gamma}>\frac{1}{2}$. Because $R_1=D_1-E_1-C_1'$ consists of $(-2)$-curves, the contribution from each other twig is at least $\frac{1}{2}$, so in fact $D_1'$ has no other twigs. In particular, $\wt Q_1$ has at most one branching component. If $\wt Q_1$ is a chain then, because $\tau_1^*=0$, $\wt Q_1=[(2)_{t_1},1]$ for some $t_1\geq 0$ and we have a contradiction as before. Thus $\wt Q_1$ is a fork and $A_0$ meets its two tips different than $C_1$. Let $B$ be the branching component. Because $C_1$ is a tip of $\wt Q_1$, the twig $T$ containing $C_1$ consists of $C_1$ and some number of $(-2)$-curves. After the contraction of $T$, $B$ becomes a $(-1)$-curve, so we have $B^2=-2$ and because $n_{exc}=0$, $B$ is not touched by $\psi$. Let $U$ be a component of $\wt Q_1-T$ meeting $B$. It is not contracted by $\psi$, so $\psi(U)$ is a $(-2)$-curve and $\psi(\wt Q_1-U-C_1+A_0)$ is a $(-2)$-chain disjoint from $E_1$. Thus the contraction of $\wt Q_1+A_0-U$ maps $X_1$ onto $\PP^2$ and $E_1$ onto a unicuspidal curve with a semi-ordinary cusp. This is a contradiction by \ref{lem:E_with_mult<=3_is_rectifiable}. Thus we may, and shall, assume further that $n=n_{exc}$.

\bcl $\#\cal C_+=0$. \ecl

Suppose $C_1$ is touched by $\psi$. By \ref{prop:(Xi,Di)_properties_for_CN}(ii) we have: $p_2(\PP^2,\bar E)=\zeta+1$, $\cal C_+=C_1'$, $K_n\cdot C_1'+\tau_1^*=0$, $\tau_j^*=0$ for $j>1$ and $R_n$ consists of $(-2)$-curves. By \ref{lem:properties_of_(Xi,Di)_for_CN}(vii) if $E_n\cdot C_1'>E_0\cdot C_1$ then some $\psi_i$ contracts $A_{i-1}$ and a $(-2)$-twig meeting $C_1$, and hence touches the image of $C_1$ once. But such $\psi_i$ would create no $(-1)$-curve in $D_i-E_i$. Therefore $$E_n\cdot C_1'=E_0\cdot C_1=\tau_1=\tau_1^*+1+s_1.$$ Blow up $m\geq 0$ times on the intersection of the proper transforms of $E_n$ and $C_1'$, denote the resulting morphism by $p\:\wt X\to X_n$ and the proper transforms of $C_1'$ and $E_n$ by $\wt C_1$ and $\wt E$ respectively. Then $(\wt C_1)^2=(C_1')^2-m$ and $$\wt C_1\cdot \wt E=\tau_1-m=3+s_1+(C_1')^2-m.$$ Take $m=(C_1')^2$. Then $\wt C_1$ is a $0$-curve, so \ref{lem:cuspidal_of_gt_and_khalf}(iv) gives $s_1=1$. Let $U$ be the component of $D_0$ which is made into a $(-1)$-curve by the (unique) $\psi_{i_0}$ touching $C_1$. Since $s_j=1$ for each $j$, $U\cdot E_0=0$. If $U$ meets $C_1$ then $\psi(U)$ is a $(-1)$-curve meeting $C_j'$ twice, so the divisor $\psi(U)+C_1'$ is nef. But in the latter case $$0\leq (\psi(U)+C_1')\cdot (2K_n+E_n)=\psi(U)\cdot E_n-\tau^*\leq 1-\tau_1^*,$$ which is impossible. Therefore $U\cdot C_1=0$. Take $m=(C_1')^2+1$. Then $\wt C_1$ is a $(-1)$-curve and it meets transversally $\wt U$, the proper transform of $U$, for which $\wt U^2=\psi(U)^2=-1$. Because $\wt U+\wt C_1$ is nef, by \ref{lem:cuspidal_of_gt_and_khalf}(iv) we get $\wt E\cdot (\wt U+\wt C_1)\geq 4$, so $E_n\cdot \psi(U)=\wt E\cdot \wt U\geq 1$. Because each $\psi_i$, $i\geq 1$ creates a $(-1)$-curve in $D_i-E_i$, \ref{lem:properties_of_(Xi,Di)_for_CN}(vii) shows that no $\psi_i$ with $i\leq i_0$ increases the intersection of images of $E_0$ and $U$, hence the image of $U$ on $X_{i_0}$ is a $(-1)$-curve disjoint from $E_{i_0}$. By \ref{lem:properties_of_(Xi,Di)_for_CN}(v) no $\psi_i$ with $i>i_0$ touches $U$, so $\psi(U)\cdot E_n=0$; a contradiction.

\med We obtain $\#\cal C_+=0$ and \beq\label{eq:claim} \tau^*+K_n\cdot R_n=p_2(\PP^2,\bar E)-\zeta\leq 1.\eeq Then $D_n-E_n$ consists of at most one $(-3)$-curve, some number of $(-2)$-curves and $c+n$ $(-1)$-curves. Moreover, by \ref{lem:properties_of_(Xi,Di)_for_CN}(vii) the intersection with the image of $E_0$ may increase under $\psi$ only for components of $D_0$ which become components of $\cal L$.

\bcl If $L$ is a component of $\cal L$ with $L\cdot E_n\leq 1$ and $M$ a component of $D_n-E_n-L$ meeting $L$ then $L\cdot M=1$. Furthermore, if $M^2=-1$ then $\tau_1^*=s_1=L\cdot E_n=1$ and the $\psi_{i_0}$ creating $L$ contracts only $A_{i_0-1}$.
\ecl

First of all note that each $\psi_i$ can increase the intersection of two components of the boundary divisor by at most one. Because $n=n_{exc}$, \ref{lem:properties_of_(Xi,Di)_for_CN}(v) implies that for any two components $V\neq W$ of $D_0$ we have $\psi(V)\cdot \psi(W)\leq V\cdot W+1$. Say $L$ is created by some $\psi_{i_0}$ for $i_0>0$. Denote the proper transforms of $M$ and $L$ on $X_0$ by $M'$ and $L'$ respectively. Put $m=-M^2$.

Assume $M$ is not a $(-1)$-curve and suppose $M\cdot L>1$. Then $M'\neq C_j$ and $M'$ meets $L'$. It follows that $M'\cdot L'=1$ and $M\cdot L=2$. By \eqref{eq:claim} $\tau^*+m\leq 3$. If $L\cdot E_n=0$ then \ref{lem:properties_of_(Xi,Di)_for_CN}(i) gives $M\cdot E_n+2m\geq 8$, so $M\cdot E_n=2$ and $m=3$, hence $\tau^*=0$. But in the latter case $M'\cdot E_0>0$, so $\tau^*>0$; a contradiction. Thus $L\cdot E_n\geq 1$. The $\psi_{i_0}$ creating $L$ is the only $\psi_i$ touching the proper transform of $L$. Indeed, this is a consequence of \ref{lem:properties_of_(Xi,Di)_for_CN}(iii),(v) and the equality $n=n_{exc}$. Since $\psi_{i_0}$ does not touch $E_{i_0-1}$, we have $L'\cdot E_0=L\cdot E_n=1$, so $\tau^*=1$ and, since $M'$ and $L'$ are contained in the same $\wt Q_j$, $M'\cdot E_0=0$. Then $M\cdot E_n\leq 1$ and $m=2$, which contradicts \ref{lem:properties_of_(Xi,Di)_for_CN}(i).

Assume $M^2=-1$. By \ref{lem:properties_of_(Xi,Di)_for_CN}(iii) $M$ is not a component of $\cal L$, so $M'=C_j$ for some $j$, say for $j=1$. Then $M+L=C_1'+L$ is nef, so $$0\leq (C_1'+L)\cdot (2K_n+E_n)=\tau_1^*+s_1+L\cdot E_n-3,$$ so, since by assumption $L\cdot E_n\leq 1$, we have $\tau_1^*=s_1=L\cdot E_n=1$. By \eqref{eq:claim} $R_n$ consists of $(-2)$-curves and $\tau_j^*=0$ for $j\geq 2$. Because $s_1=1$, $L'$ is disjoint from $E_0$ and $A_{i_0-1}$ meets $E_{i_0-1}$. Since $C_1$ is not touched by $\psi$, $L'$ is the unique component of $D_0-E_0$ meeting $C_1$, hence a $(-2)$-curve. If $A_{i_0-1}$ meets $L'$ then $\psi_{i_0}$ contracts only $A_{i_0-1}$ and we are done. Thus we may assume $A_{i_0-1}$ does not meet $L'$. Then there is a $(-2)$-twig $\Delta_L$ in $D_0-E_0$ meeting $L'$ contracted by $\psi$. Since $C_1$ is not touched by $\psi$, $L'$ is a branching $(-2)$-curve in $D_0-E_0$. Since $K_n\cdot R_n=0$ and $n=n_{exc}$, \ref{lem:properties_of_(Xi,Di)_for_CN}(vii) implies that images of non-contracted components of $\wt Q_1-C_1-L'$ are either $(-1)$-curves meeting $E_n$ at most once or $(-2)$-curves disjoint from $E_n$. The divisor $\wt Q_1$ is a rational snc-tree contractible to a point, so $L'$ meets a component $U$ with self-intersection smaller than $(-2)$ (cf. \ref{lem:Q_with_small_KQ}). Because $L'$ has self-intersection $(-2)$, $\psi$ touches it once, so $\psi_*(U)\neq 0$. By \ref{prop:(Xi,Di)_properties_for_CN}(iii) $\psi(U)^2\neq -1$, so $\psi(U)^2=-2$. But $U^2<-2$, so $\psi$ touches $U$. Since $n=n_{exc}$, there is a component $V$ of $\cal L-L$ meeting $\psi(U)$. By \ref{lem:properties_of_(Xi,Di)_for_CN}(vii) for the nef divisor $f=L+\psi(U)+V$ we obtain $f\cdot E_n=1+V\cdot E_n\leq 2+\psi_*^{-1}(V)\cdot E_0\leq 3$.  Then $f\cdot (2K_n+E_n)<0$; a contradiction.

\bcl $\zeta=p_2(\PP^2,\bar E)-1$. \ecl

Suppose $\zeta=p_2(\PP^2,\bar E)$. By \eqref{eq:claim} $R_n$ consists of $(-2)$-curves and $\tau^*=0$. Then $C_j$'s are the only components of $D_0$ meeting $E_0$ and they are not touched by $\psi$. Suppose $L$ is a $(-1)$-curve in $D_n$ created by $\psi$. By definition, it is not a superfluous component of $D_n$. By \ref{lem:properties_of_(Xi,Di)_for_CN}(vii) $L\cdot E_n\leq 1$, so the previous claim implies that we can find two $(-2)$-curves, say $V_1$, $V_2$ in $D_n$ meeting $L$. Moreover, by \ref{lem:properties_of_(Xi,Di)_for_CN}(vii) $V_1, V_2$ are disjoint from $E_n$, so for $f=V_1+2L+V_2$ we have $f\cdot (2K_n+E_n)=2L\cdot E_n-4<0$. But $2K_n+E_n\geq 0$ and $f$ is nef; a contradiction. Thus $n_{exc}=n=0$ and $R_0$ consists of $(-2)$-curves. It follows that $D_0-E_0$ contains no branching components, so all cusps of $\bar E$ are semi-ordinary. Such curves have been ruled out in \ref{lem:E_with_mult<=3_is_rectifiable}; a contradiction.

\bcl $K_0\cdot R_0=K_n\cdot R_n=1$. \ecl

We first show that $K_0\cdot R_0=K_n\cdot R_n$. Because every $\psi_{i+1}$, $i\geq 0$ creates a $(-1)$-curve in $D_i-E_i$, we may assume that there exists $i\geq 0$, such that $A_i\cdot E_i=0$, otherwise we are done. Let $L\leq \cal L$ be the image of the $(-1)$-curve created by $\psi_{i+1}$. Then $L\cdot E_n=L'\cdot E_0\leq 1$, where $L'=\psi_*^{-1}(L)$. By the third claim we may also assume $L$ does not meet $(-1)$-curves in $D_n-E_n-L$. Then, as before, we can find in $D_n-E_n$ some $(-2)$-curves $V_1$ and $V_2$ meeting $L$. We have $n=n_{exc}$, so since $V_i$ are not $(-1)$-curves, by \ref{lem:properties_of_(Xi,Di)_for_CN}(vii) $(V_1+V_2)\cdot E_n=\psi_*^{-1}(V_1+V_2)\cdot E_0$. Because $s_j=1$ for $j\neq 1$, we have $\psi_*^{-1}(V_1+V_2)\cdot E_0\leq 1$, so for the nef divisor $f=V_1+2L+V_2$ we get $0\leq f\cdot (2K_n+E_n)\leq\psi_*^{-1}(V_1+V_2)\cdot E_0-2<0$; a contradiction. Thus $K_0\cdot R_0=K_n\cdot R_n\leq 1$. If $K_0\cdot R_0=0$ then $D_0-E_0$ has no branching components and hence the cusps $q_2,\ldots,q_c$ are semi-ordinary and $q_1$ has multiplicity three, which is impossible by \ref{lem:E_with_mult<=3_is_rectifiable}.

\med Now the equation \eqref{eq:claim} gives $\tau^*=0$. Because we know that each $\wt Q_j$ contracts to a smooth point and that $C_j$ is a tip of $\wt Q_j$, it follows that the cusps $q_2, q_3,\ldots,q_c$ are semi-ordinary and $\wt Q_1$ is a fork with a $(-2)$-curve as a branching component and maximal twigs of types $[2], [(2)_{t_1},3]$ and $[1,(2)_{t_2}]$ for some $t_1,t_2\geq 0$. Let $U$ be the $(-3)$-curve in $\wt Q_1$. The argument with finding nef chains of type $[1,2,\ldots,2,1]$ in $D_n-E_n$ shows that no $A_{i}$ meets $\wt Q_2+\ldots+\wt Q_c$ or the twigs $[(2)_{t_2}]$ or $[2]$. If $n\neq 0$ then, because $K_n\cdot R_n=1$, $\psi$ does not touch $U$, so $A_0$ meets the component of $[(2)_{t_1}]$ meeting $U$ and $t_1\geq 2$. But in the latter case the contraction of $D_0-E_0-U+A_0$ maps $X_0$ to $\PP^2$ and $E_0$ to a cuspidal curve with semi-ordinary cusps, which contradicts \ref{lem:E_with_mult<=3_is_rectifiable}.

Thus $n=0$. We obtain $\#D_0=\rho(X_0)=10-K_0^2=10+K_0\cdot E_0-\zeta,$ so $\#D_0=9+\gamma_0-p_2(\PP^2,\bar E)$ and hence $x+t_1+t_2=4+\gamma_0-p_2(\PP^2,\bar E)$, where $x=\#(\wt Q_2+\ldots+\wt Q_c)$. In particular, $t_1\leq 4+\gamma_0-p_2(\PP^2,\bar E)$. The sequence of characteristic pairs for $\wt Q_1$ is $\binom{2t_1+3}{2}, \binom{1}{1}_{t_2+1}$, so the equations in \ref{cor:I and II}(i),(ii) read as \begin{eqnarray*} 3d &=&2t_1-2p_2(\PP^2,\bar E)+\gamma_0+20\\ d^2 &=&12t_1+44-4p_2(\PP^2,\bar E)+3\gamma_0.\end{eqnarray*} Now \ref{prop:(Xi,Di)_properties_for_CN}(i) gives $4\leq \gamma_0\leq 4p_2(\PP^2,\bar E)-2$ and one checks that with this bound solutions exist only for $p_2(\PP^2,\bar E)=4$. By \ref{lem:(Xi,Di)_properties}(iv) $\ind(D)\leq 1$. But $Q_1$ has three maximal twigs and at least two of them are $(-2)$-twigs, so the contribution to $\ind(D)$ from $Q_1$ is bigger than~$1$; a contradiction.
\end{proof}

\blem \cite[Theorem 3.3]{Kumar-Murthy}. \label{lem:elliptic_ruling} Let $E$ be a smooth rational curve on a smooth rational surface $Y$. If $E^2=-4$ and $C$ is a $(-1)$-curve for which $E\cdot C=2$ then $|E+2C|$ induces an elliptic fibration of $Y$. In particular, if $2K_Y+E\sim 0$ then any $(-1)$-curve on $Y$ gives such a fibration. Moreover, in the latter case the fibration has no section and singular fibers other than $E+2C$ consist of $(-2)$-curves. \elem

\begin{proof} The divisor $f=E+2C$ has arithmetic genus $1$ and self-intersection zero. We have $E+2C=(E+C+K_Y)+(C-K_Y)$ and $h^0(E+C+K_Y)\geq 1$, $h^0(C-K_Y)\geq 1$ by Riemann-Roch. Since $C$ is not in the fixed part of $E+C+K_Y$ or $C-K_Y$, there is more than one effective divisor in $|E+2C|$, so $h^0(E+2C)\geq 2$. The exact sequence $$0\to \cal O_{Y}(2C)\to \cal O_{Y}(E+2C)\to \cal O_{E}\to 0$$ gives in fact $h^0(E+2C)=2$. This gives an elliptic fibration $\theta\:Y\to\PP^1$ with $E+2C$ as a fiber. If the equivalence $2K_Y+E\sim 0$ holds then $\theta$ has no $1$-section because $f\sim 2(C-K_Y)$ and a component of a singular fiber other than $E+2C$ intersects $2K_Y\sim -E$ trivially, hence it is a $(-2)$-curve.
\end{proof}

\blem\label{lem:h0=2_only_for_zeta=1} If $p_2(\PP^2,\bar E)=2$ then $\zeta=1$. \elem

\begin{proof} Suppose $\zeta=K_n\cdot (K_n+E_n)\neq 1$. By \ref{prop:jump>=2} $\zeta\leq 0$, so \ref{prop:(Xi,Di)_properties_for_CN}(i) gives $4\leq  \gamma_n+\tau^*+(n-\eta)\leq 2p_2(\PP^2,\bar D)+2\zeta\leq 4$. Thus $\zeta=0$, $\gamma_n=4$, $\tau^*=0$ and $n=\eta$. In particular, $n=n_{exc}$ and by \ref{lem:properties_of_(Xi,Di)_for_CN}(vi) $\#\cal C_+=0$. Then \ref{prop:(Xi,Di)_properties_for_CN}(ii) gives $K_n\cdot R_n=2$. Note that by \ref{lem:properties_of_(Xi,Di)_for_CN}(vii) for every component $V$ of $R_n$ we have $V\cdot E_n=\psi_*^{-1}(V)\cdot E_0=0$. Even with all this information the proof is long.

\bcl Every component of $\cal L$ meets $E_n$ once. \ecl

Suppose $L\cdot E_n\neq 1$ for some component $L$ of $\cal L$. Put $L'=\psi_*^{-1}(L)$. By \ref{lem:(Xi,Di)_properties}(vi) $L\cdot E_n\leq L'\cdot E_0+1=1$, so $L\cdot E_n=0$. Suppose first that every component of $D_n-E_n$ meets $L$ at most once. Since $L$ is not superfluous, at least three components of $D_n-E_n$, say $V_1, V_2, V_3$, meet $L$. By \ref{lem:properties_of_(Xi,Di)_for_CN}(iii) they are components of $R_n+\cal C_{exc}+\cal C_+$. However, if some $C_j'$, say $C_1'$, meets $L$ then $f=L+C_1'$ is nef and $f\cdot(2K_n+E_n)=C_1'\cdot E_n-4=C_1\cdot E_0-4=\tau_1-4<0$, which is impossible, because $2K_n+E_n\geq 0$. We infer that $V_i$ are components of $R_n$, so they do not meet $E_n$. Only one of them may be a $(-2)$-curve. Indeed, if $V_1$ and $V_2$ are $(-2)$-curves then $f=V_1+2L+V_2$ is nef and $0\leq f\cdot (2K_n+E_n)=-4$; a contradiction. Since $K_n\cdot R_n=2$, \ref{lem:properties_of_(Xi,Di)_for_CN}(v) implies that $V_1^2=V_2^2=-3$ and $V_3^2=-2$. The divisor $f=3L+V_1+V_2+V_3$ is nef. Because $2L$ is in the fixed part of $2K_n+E_n$, we have $0\leq f\cdot (2K_n+E_n-2L)=2(V_1+V_2+V_3)\cdot(K_n-L)=-2$; a contradiction.

Therefore, there is a component $V'\leq D_0$, say $V'\leq \wt Q_1$, such that for $V=\psi(V')$ we have $V\cdot L\geq 2$. Again, $V$ is a component of $R_n$ and hence is disjoint from $E_n$. As in the proof of the third claim of \ref{prop:jump>=2} we note that $V\cdot L\leq V'\cdot L'+1$, so $L'\cdot V'=1$ and $L\cdot V=2$. By \ref{lem:properties_of_(Xi,Di)_for_CN}(i) $V^2=-4$, and hence $R_n-V$ consists of $(-2)$-curves.

Suppose some other $(-1)$-curve $U\leq D_n-E_n$ meets $V$. Then $f=V+2L+U$ is nef, hence $0\leq f\cdot (2K_n+E_n)=U\cdot E_n-2$, so $U=C_1$. Because $C_1$ is a tip of $\wt Q_1$ not touched by $\psi$, $C_1\cdot V'\neq 0$, and then $V'$ is a $(-2)$-curve. However, $(V')^2\leq V^2\leq -4$; a contradiction. Thus $L$ is the unique $(-1)$-curve in $D_n$ meeting $V$. Therefore, there is a component $M$ of $R_n$ meeting $L+V$. In fact $M\cdot L=0$, otherwise $f=V+3L+M$ is a nef divisor intersecting $2K_n+E_n$ negatively. Because $n=n_{exc}$, the fact that $L$ is the unique $(-1)$-curve in $D_n-E_n$ meeting $V$ implies that $V\cdot M=V'\cdot \psi_*^{-1}(M)\leq 1$, so $V\cdot M=1$. Put $f_1=V+2L$ and $f_2=E_n+2C_1'$. We have $f_1\cdot f_2=2C_1'\cdot V+4C_1'\cdot L=0$. By \ref{lem:elliptic_ruling} $f_1$ and $f_2$ are fibers of the same elliptic fibration of $X_n$, so $f_1\sim f_2$. But $f_1\cdot M=1$ and $f_2\cdot M=2C_1'\cdot M$, which is even; a contradiction.

\bcl $K_0\cdot R_0=2$ and each $\psi_i$ with $i\geq 1$ contracts exactly $A_{i-1}$. \ecl

We may assume $n\neq 0$. Let $L$ be any component of $\cal L$. Let $\psi_{i_0}$ be the first $\psi_i$ with $i\geq 1$ touching some proper transform of $L$. In fact it is $\psi_{i_0}$ which creates $L$. Indeed, because $n=n_{exc}$, $\psi_{i_0}$ creates a $(-1)$-curve, whose image on $X_n$ would otherwise be a $(-1)$-curve (cf. \ref{lem:properties_of_(Xi,Di)_for_CN}(v)) meeting $L$, contradicting \ref{lem:properties_of_(Xi,Di)_for_CN}(iii). It follows that only $\psi_{i_0}$ touches some proper transform of $L$. Put $L'=\psi_*^{-1}(L)$. We have $L'\cdot E_0=0$ and $L\cdot E_n=1$, so by \ref{lem:properties_of_(Xi,Di)_for_CN}(vii) $A_{i_0-1}\cdot E_{i_0-1}=1$. Then $\psi_{i_0}$ is of type $I$, so $\psi_{i_0}$, and hence $\psi$, touches $L'$ once, so $L'$ is a $(-2)$-curve. It follows that $K_0\cdot R_0=K_n\cdot R_n=2$.

Suppose some $\psi_i$ with $i\geq 1$ contracts more then just $A_{i-1}$. Then $\psi_i$ contracts a maximal $(-2)$-twig of $D_{i-1}-E_{i-1}$ and hence touches a component of $D_{i-1}-E_{i-1}$, say $B$, meeting it. Because $K_0\cdot R_0=K_n\cdot R_n$, $B$ is a branching $(-2)$-curve in $D_{i-1}-E_{i-1}$. Now the equality $n=\eta$ implies that $B$ meets another $(-2)$-twig of $D_{i-1}-E_{i-1}$, which contradicts \ref{lem:(Xi,Di)_properties}(v).

\med We have $K_n^2=\zeta-K_n\cdot E_n=-2$, so the Noether formula gives $\rho(X_n)=10-K_n^2=12$ and hence $\#D_n=12+n$ by \ref{lem:(Xi,Di)_properties}(iii). Since every $A_i$ meets $E_i$, $\gamma_0=\gamma_n+n$. The claim above gives $$\#D_0=\#D_n=12+n.$$ We have $K_0\cdot R_0=2$. Let $U_1$ be a component of $R_0$ with $U_1^2\leq -3$. If $U_1^2=-4$ then such $U_1$ is unique (then we put $U_2=0$). If $U_1^2=-3$ then there is another such curve, call it $U_2$, and $U_2^2=-3$. The divisor $R_0-U_1-U_2$ consists of $(-2)$-curves. The morphism $\psi$ does not touch $U_1+U_2$. Note that if $n\geq 2$ then, since $\eta=n$, $D_n$ has at least two and hence $D_n'$ has at least three $(-2)$-twigs, so $\ind(D_n')\geq \frac{3}{2}$. But \ref{lem:(Xi,Di)_properties}(iv) gives $n+\ind(D_n')\leq 3$, so $n\leq 1$. Observe also that if $f=[1,2,\ldots,2]$ is any chain in $D_0-E_0$ containing come $C_j$ then $\psi$ does not touch it. Indeed, otherwise, because $n=n_{exc}$, $D_n-E_n$ contains a nef chain of type $[1,2,\ldots,2,1]$ and such a chain meets $E_n$ at most three times, which contradicts \ref{lem:cuspidal_of_gt_and_khalf}(iii). In particular, the cusps $q_j$ for which $\wt Q_j$ contains neither $U_1$ nor $U_2$ are semi-ordinary, so in \ref{cor:I and II} applied to $(X_0,D_0)$ we obtain $M(q_j)=I(q_j)=\#\wt Q_j$. Let $x$ be the total number of components of all such $\wt Q_j$'s. The equations in \ref{cor:I and II}(i),(ii) read as \begin{eqnarray}\label{eq:pairs-case_(p2,zeta)=(2,0),I} 3d'+1&=&\sum_j M(q_j)\\ \label{eq:pairs-case_(p2,zeta)=(2,0),II}(d')^2+1&=&\sum_j I(q_j), \end{eqnarray} where $d'=d/2$.

\bcl $U_2\neq 0$. $U_1$ and $U_2$ belong to different connected components of $D_0-E_0$. \ecl

We may assume $U_1$ is a component of $\wt Q_1$. Suppose $U_2=0$. Then $U_1^2=-4$ and $\wt Q_1$ is a fork with maximal twigs $[2,2]$, $[(2)_{t_1},4]$ and $[1,(2)_{t_2}]$ for some $t_1,t_2\geq 0$ (cf. \ref{lem:Q_with_small_KQ}(iii)). The characteristic pairs of $\wt Q_1$ are $\binom{3t_1+4}{3}$, $\binom{1}{1}_{t_2+1}$. We have $12+n=\#D_0=t_1+t_2+x+6$, so $3d'=2t_1+n+12$ and $(d')^2=8t_1+n+18$, hence $(d')^2-12d'+30+3n=0$. The latter equation is equivalent to $(d'-6)^2=3(2-n)$, which has no solutions for $n\leq 1$; a contradiction.

Thus $U_2\neq 0$, so $U_1^2=U_2^2=-3$. Let $T$ be the maximal twig of $\wt Q_1$ containing the $(-1)$-tip $C_1$. Since $\wt Q_1$ can be contracted to a smooth point by successive contractions of $(-1)$-curves, we see that $T-C_1$ consists of $(-2)$-curves and meets $\wt Q_1-T$ in a branching $(-2)$-curve. Suppose $\wt Q_1-T$ is a chain. Then $\wt Q_1-T=[3,1,2,3,(2)_{t_1}]$ for some $t_1\geq 0$. The sequence of characteristic pairs for $\wt Q_1$ is $\binom{3t_1+5}{3},\binom{1}{1}_{t_2+1}$ for some $t_2\geq 0$ and $12+n=\#D_0=t_1+t_2+x+6$, so \eqref{eq:pairs-case_(p2,zeta)=(2,0),I} and \eqref{eq:pairs-case_(p2,zeta)=(2,0),II} combine to give $(d')^2-3d'=6t_1+8$. The latter equation gives $3|(d')^2+1$; a contradiction.

It follows that $\wt Q_1-T$ is a fork with a branching component $B$ which is a $(-2)$- or a $(-3)$-curve. In the first case the maximal twigs of this fork are $[2]$, $[(2)_{t_1},3]$, $[2,2,3,(2)_{t_2-1}]$ for some $t_1\geq0$, $t_2\geq 1$ and in the second case $[2]$, $[(2)_{t_1},3]$, $[2,2]$. The sequence of characteristic pairs for $\wt Q_1$ is $$\binom{4t_1+6}{4},\binom{2}{2}_{t_2},\binom{2}{1},\binom{1}{1}_{t_3+1}$$ for some $t_1,t_2,t_3\geq 0$. Also, $\#D_0=t_1+t_2+t_3+x+7$, so $t_1+t_2+t_3+x=n+5$. The equations \eqref{eq:pairs-case_(p2,zeta)=(2,0),I} and \eqref{eq:pairs-case_(p2,zeta)=(2,0),II} read as \begin{eqnarray*}3d'&=&3t_1+t_2+16+n\\ (d')^2&=&15+3t_2+34+n.\end{eqnarray*} But $t_1+t_2\leq n+5\leq 6$ and with this bound they have no solutions; a contradiction. Thus, $U_2$ belongs to, say, $\wt Q_2$.

\bcl $n=0$ and $\deg \bar E\in \{8,10\}$. \ecl

For $j=1,2$ the curve $C_j$ is a tip of $\wt Q_j$ and $\wt Q_j-C_j-U_j$ consists of $(-2)$-curves. Then $\wt Q_j$ is a fork with a branching $(-2)$-curve and maximal twigs of type $[2],[(2)_{t_j},3], [1,(2)_{t_j'}]$ for some $t_j, t_j'\geq 0$. We obtain $\#\wt Q_1+\#\wt Q_2+1+x=\#D_0=12+n$, so $t_1+t_1'+t_2+t_2'+x=n+3\leq 4$. The sequence of characteristic pairs of $\wt Q_j$ for $j=1,2$ is $\binom{2t_j+3}{2}, \binom{1}{1}_{t_j'+1}$, so $3d'=t_1+t_2+n+12$ and $(d')^2=3(t_1+t_2)+n+16$, where $d'=\deg \bar E/2$. In particular,the second equation taken modulo $3$ gives $n=0$, so $t_1+t_2+t_1'+t_2'+x=3$. The solutions are $(\deg \bar E,t_1+t_2)=(8,0)$ and $(10,3)$. This is one of the most 'resistant' cases in the article.

\bcl There is a $(-1)$-curve $V$ for which $V\cdot (D_0-U_1-U_2)=0$ and $V\cdot U_1=V\cdot U_2=1$. \ecl

We have $2K_0+E_0\geq 0$ and $(2K_0+E_0)^2=4(\zeta -1)<0$, so there is an irreducible curve $V$ such that $2K_0+E_0-V\geq 0$ and $(2K_0+E_0)\cdot V<0$. Clearly, $V\neq E_0$, so we get $V^2<0$ and $K\cdot V<0$, hence $V$ is a $(-1)$-curve meeting $E_0$ at most once. But $E_0\cdot (2K_0+E_0)=\gamma_0-4=0$, so since $E_0$ is not in the fixed part of $2K_0+E_0-V$, we have in fact $V\cdot E_0=0$. Since the inequality \ref{prop:(Xi,Di)_properties_for_CN}(i), which is an equality in our case, is equivalent to $(2K_0+D_0^\flat)\cdot (2K_0+E_0)\geq 0$ ($n=0$, so the divisor $2K_0+D_0^\flat$ is nef), we get $(2K_0+D_0^\flat)\cdot V=0$. Therefore $$V\cdot (D_0-E_0-\Upsilon_0-\Delta^+_0-\Bk'\Delta_0^-)=2.$$ The argument with finding nef divisors in $D_0-E_0$ of type $[1,2,\ldots,2,1]$ shows that $V$ is disjoint from $D-\Delta_0^--U_1-U_2$ and meets $\Delta_0^-$ at most once. We get the equality $V\cdot (D_0-E_0-\Bk'\Delta_0^-)=2$. In particular, $V\cdot \Bk'\Delta_0^-$ is an integer, so $V$ is disjoint from $\Delta_0^-$, hence $V\cdot D=V\cdot(U_1+U_2)=2$. If $V$ meets only one $U_i$ then $U_i+2V$ is nef and intersects $2K_0+E_0$ negatively, which is impossible. We obtain $V\cdot U_1=V\cdot U_2=1$.

\med Let $p\:(X_0,D_0)\to (Y,D_Y)$ be the contraction of $V$. The dual graph of $p_*(\wt Q_1+\wt Q_2)$ is $$ \xymatrix{
{-1}\ar@{-}[r] &{(-2)_{t_1'}}\ar@{-}[r] &{-2}\ar@{-}[r]\ar@{-}[d] &{-2}\ar@{-}[r]\ar@{-}[d]  &{-2}\ar@{-}[r]\ar@{-}[d]  &{-2}\ar@{-}[r]\ar@{-}[d] &{(-2)_{t_2'}}\ar@{-}[r] &{-1}
\\ {} & {} & {-2} & {\ \ (-2)_{t_1}} & {\ \ (-2)_{t_2}} & {-2}& {} & {}}$$
where $t_1+t_2+t_1'+t_2'+x=3$ and $t_1+t_2\in\{0,3\}$.

\bcl $x=t_1=t_2=0$ and $t_1'+t_2'=3$. \ecl

Because $\psi_0$ does not touch $V$, the proper transform $V'$ of $V$ on $X$ is a $(-1)$-curve. Let $p'\:(X,D)\to (Y',D_{Y'})$ be the contraction of $V'$. Note $D_Y'$ is an snc-divisor. The surface $Y'\setminus D_{Y'}$ contains no $\C^1$'s, because $X\setminus D$ contains none. Then \ref{lem:Bk} gives $((K_{Y'}+D_{Y'})^-)^2=-\ind (D_{Y'})$. We have $(p')^*(K_{Y'}+D_{Y'})=K+D+V$, so by \ref{lem:(Xi,Di)_properties}(ii) $(K_{Y'}+D_{Y'})^2=(K+D)^2+(K+D)\cdot V=p_2(\PP^2,\bar E)-2+1=1,$ giving $((K_{Y'}+D_{Y'})^+)^2=\ind (D_{Y'})+1$. We compute $\chi(Y'\setminus D_{Y'})=\chi(X\setminus D)-\chi(V\cap (X\setminus D))=1$. Then the logarithmic Bogomolov-Miyaoka-Yau inequality for $(Y',D_{Y'})$ reads as $\ind (D_{Y'})\leq 2$. But $D_{Y'}$ is obtained from $D_Y$ by resolving the tangency of $C_j$'s and $E_0$, which replaces the chains $[(2)_{t_j'}]$ with $[2,1,3,(2)_{t_j'}]$. The contribution to $\ind (D_{Y'})$ from the twigs of $D_{Y'}$ over $q_j$, $j=1,2$ is therefore $2\cdot\frac{1}{2}+\frac{t_j}{t_j+1}$, so the inequality $\ind (D_{Y'})\leq 2$ shows that $t_1=t_2=0$ and there is no contribution from $\wt Q_3+\ldots+\wt Q_c$. The latter means that $x=0$, which gives $t_1'+t_2'=3$.

\med By symmetry we may assume $t_1'\in\{2,3\}$. Let $G_1$ be the component of $p_*\wt Q_1$ meeting $p_*C_1$ and let $\mu F$, where $\mu$ is the multiplicity, be the fiber of the elliptic fibration of $Y$ induced by $|p_*(E_0+2C_1)|$ (see \ref{lem:elliptic_ruling}) which contains $p_*(\wt Q_1+\wt Q_2-C_1-C_2)-G_1$. Since $K_Y^2=K_0^2+1=-1$, all fibers other than $p_*(E_0+2C_1)$ are minimal, i.e.\ they consist of $(-2)$-curves. Suppose $t_1'=2$. Then the latter divisor has two branching components, hence equals $F_{red}$, and $F$ is of Kodaira type $I^*_3$ (called also $\wt D_7$; for the classification of fibers and basic facts on elliptic fibrations see for example \cite[V.7]{BHPV}). Because $F_{red}$ is simply connected, $\mu=1$. Since $G_1$ meets one of the tips of $F_{red}$, whose multiplicities in $F$ are $1$, we get $G\cdot (\mu F)=G\cdot F_{red}=1$. This is a contradiction, because $G_1$ is a $2$-section of the fibration. Therefore $t_1'=3$. In this case $p_*(\wt Q_1+\wt Q_2-C_1-C_2)-G_1$ is a $(-2)$-fork with maximal twigs of lengths $1,2,4$, which implies that $F$ is of Kodaira type $\textrm{II}^*$ ($\wt E_8$). But then the component of $F$ meeting $C_2$ has multiplicity three in $F$. This is in contradiction with the fact that $C_2$ is a $2$-section.
\end{proof}

We now prove Theorem \ref{thm:MAIN1}.

\begin{proof}[Proof of Theorem \ref{thm:MAIN1}] By \ref{prop:jump>=2} and \ref{prop:(Xi,Di)_properties_for_CN}(i) $4p_2(\PP^2,\bar E)-2\geq 2\zeta+2p_2(\PP^2,\bar E)\geq \gamma_n\geq 4,$ so $p_2(\PP^2,\bar E)\geq 2$. By \ref{lem:(Xi,Di)_properties}(iv) $p_2(\PP^2,\bar E)\leq 4$. Because $(K+D)^2=K\cdot(K+D)-2=p_2(\PP^2,\bar E)-2$, we only need to show that $p_2(\PP^2,\bar E)\neq 2$. Suppose $p_2(\PP^2,\bar E)=2$.

\bcl $n=\zeta=1$, $n_{exc}=n_1=0$ and $2K_1+E_1\sim 0$. \ecl

Lemma \ref{lem:h0=2_only_for_zeta=1} gives $\zeta=1$, so by \ref{prop:jump>=2} $n=1$, $n_{exc}=n_1=0$ and $2K_1+E_1\sim 0$.

\med Note that in particular $\eta=0$. We have also $K_1\cdot E_1=2\zeta=2$, so $\gamma_1=4$ and $K_1^2=-1$. By \ref{prop:(Xi,Di)_properties_for_CN}(ii) $\#\cal C_+=\tau^*=K_1\cdot R_1=0$. By \ref{lem:(Xi,Di)_properties}(iii) $\#D_1=\rho(X_2)+1=11-K_1^2=12$. We may assume $A_0$ meets $\wt Q_1$ and does not meet $\wt Q_3+\ldots+\wt Q_c$. Then the cusps $q_3,\ldots, q_c$ are semi-ordinary.

\bcl $A_0$ does not meet $\wt Q_j$ for $j\geq 2$.\ecl

Suppose $A_0$ meets $\wt Q_2$. Let $j=1,2$. Because $\psi=\psi_1$ does not create $(-1)$-curves in $D_1-E_1$, $\wt Q_j$ contains a component with self-intersection smaller than $(-2)$ and $\psi$ contracts some component of $\Delta_0^-$. In particular, $\wt Q_j$ contains a branching component. Let $T_j$ be the twig of $\wt Q_j$ (not necessarily maximal, possibly empty) contracted by $\psi$ and let $V_j\leq \wt Q_j-T_j-C_j$ be the component of $\wt Q_j$ meeting $T_j+A_0$. We see that $V_j^2\leq -3$ and that $\wt Q_j-C_j-V_j-T_j$ consists of $(-2)$-curves and is not touched by $\psi$. Contract successively $(-1)$-curves in $\wt Q_j$ until $V_j$ becomes a $(-1)$-curve; denote the images of $\wt Q_j$ and $V_j$ by $Q_j'$ and $V_j'$ respectively. Put $W_j=Q_j'-V_j'-T_j$. Clearly, $W_j$ consists of $(-2)$-curves and meets the $(-1)$-curve $V_j'$. It follows that $W_j$ is a chain; write $W_j=[(2)_{w_j}]$ for some $w_j\geq 0$. Now \ref{lem:chains_with_d=1} implies that $Q_j'=[(2)_{w_j},1]$ or $Q_j'=[(2)_{w_j},1,w_j+2,(2)_{t_j}]$ for some $t_j\geq 0$, so $T_j=0$ or $T_j=[(2)_{t_j},w_j+2]$. In particular, we see that $R_0=D_0-E_0-\cal C$ contains at most four components which are not $(-2)$-curves. Put $v_j=-V_j^2$. We have $v_j\geq 3$. The divisor $\wt Q_j-T_j$ is a fork with a branching $(-2)$-curve and maximal twigs $W_j+V_j$, $[(2)_{v_j-2}]$ and $[1,(2)_{a_j}]$ for some $a_j\geq 0$. We compute that the contribution to $\ind(D_1')$ from the maximal twigs of $D_1'$ over $q_j$ is $\frac{1}{2}+\frac{w_j}{w_j+1}+\frac{v_j-2}{v_j-1}\geq 1$. By \ref{lem:(Xi,Di)_properties}(iv) $\ind(D_1')\leq 2$. This gives $v_1=v_2=3$ and $w_1=w_2=0$. By \eqref{eq:Ai} $A_0$ meets exactly one $(-2)$-twig of $D_0$, so we have $T_1\neq 0$ and $T_2=0$. By the definition of $V_j$'s $\psi$ maps the chain $V_1+T_1+A_0+V_2$, which is of type $[3,(2)_{t_1+1},1,3]$, onto $[2,2]$ by inner contractions. But for $t_1\geq 0$ this is impossible; a contradiction.

\med The divisor $\psi_*(\wt Q_1+A_0-C_1)$ has arithmetic genus one and consists of $(-2)$-curves, so it contains a $(-2)$-cycle. Thus there exist $r\geq 2$ and components $U_1,\ldots,U_r$ in $\wt Q_1+A_0$, such that $\psi(U_i)$ are $(-2)$-curves and for $f=\psi(U_1+\ldots+U_r)$ we have $\psi(U_i)\cdot (f-\psi(U_i))=2$ for every $i$. In particular, $f^2=0$. Denote by $G_1$ the component of $\wt Q_1$ meeting $C_1$. Note that $C_1$ is a tip of $\wt Q_1$, so $G_1^2=-2$ and $\wt Q_1-G_1-C_1$ has at most two connected components, only one of which may be branched. By \ref{lem:elliptic_ruling} for each $j\leq c$ the divisor $f_j=E_1+2C_j'$ is a fiber of an elliptic fibration $\theta_j$, whose other fibers consist of $(-2)$-curves. Because $\psi$ does not touch $C_1$, it does not contract $G_1$. But because $G_1^2=\psi(G_1)^2=-2$, $\psi$ does not even touch $G_1$.

\bcl $\psi(G_1)$ is a component of $f$. \ecl

Suppose not. Then $f$ is vertical for the fibration $\theta_1$ and $\psi(G_1)$ is a $2$-section. By the Kodaira classification we infer that $2f$ is a fiber of $\theta_1$. We have $c=1$, otherwise $C_2$ would be a $2$-section of the fibration, which is impossible, because it is disjoint from $\psi(\wt Q_1)$. Furthermore, $\psi(G_1)$ is the unique horizontal component of $D_1$. Some $\psi(U_i)$, say $\psi(U_1)$, meets $\psi(G_1)$. Then $U_1$ meets $G_1$. Since $f$ is a cycle, $\psi^*(f)_{red}$ is a cycle containing $A_0$. Because the divisor $\wt Q_1-\psi^*(f)_{red}$ is not touched by $\psi$, it consists of $C_1$ and some number of $(-2)$-curves. But $G_1$ touches $U_1$, so in fact it is a chain of type $[(2)_t,1]$ for some $t\geq 1$. Let $Q_1'$ and $U_1'$ be the images of $\wt Q_1$ and $U_1$ after the contraction of the latter chain. Then $Q_1'$ is a chain and $U_1'$ is its unique $(-1)$-curve. The image of $A_0$, which is a $(-1)$-curve, meets $Q_1'$ in tips. We know that $Q_1'$ is of one of two types in \ref{lem:chains_with_d=1}. The fact that $\psi$ maps $Q_1'-U_1'$ onto a nonempty chain of $(-2)$-curves implies (an easy induction) that if $Q_1'$ is of the first type then $m_2=x$ and $m_{2i-1}=0$, $m_{2i+2}=x-1$ for $i\geq 1$ if $x\neq -1$ and $k=m_1=1$ otherwise. Similarly, for the second type we get $m_1=x+1$ and $m_{2k+1}=x-1$, $m_{2k}=0$ for $i\geq 1$ if $y\neq 0$ and $k=m_1=1$ otherwise. If, extending our previous notation, by $(a,b)_m$ we denote the sequence $(a,b,a,b,\ldots,a,b)$, where the pair is repeated $m$ times, then $Q'$ is in both cases of type $[x+3,2,1,3,(2)_x]$ or $$[(x+2)_k,x+3,1,(2)_{x+1},3,((2)_{x-1},3)_k,(2)_x]$$ for some $k\geq 0$ and $x\geq 1$. We infer that $U_1$ is not touched by $\psi$, so it is a $(-2)$-curve, which gives $t=1$. Let $U_2$ denote the component of $\wt Q_1$ meeting $U_1$ which has self-intersection $-x-3$. It is not contracted by $\psi$, so the contraction of $\wt Q_1+A_0-U_2$ maps $X_0$ onto $\PP^2$ and $E_0$ onto a unicuspidal curve for which the cusp is semi-ordinary. This contradicts \ref{lem:E_with_mult<=3_is_rectifiable}.

\bcl $\wt Q_1$ is a fork with $G_1$ as the branching component. \ecl

Because $G_1$ is not touched by $\psi$, the claim above implies that it is a branching component of $\wt Q_1$. Since $\wt Q_1$ contracts to a smooth point, there is a maximal twig of $\wt Q_1$, say $T$, meeting $G_1$. The curve $\psi(G_1)$ is a component of $f$, so $A_0$ meets $T$. Suppose there is another branching component in $\wt Q_1$. Say $V$ is the one nearest to $G_1$ inside $D_0$. Let $W_1+\ldots+W_z$, $z\geq 1$ be the shortest chain in $\wt Q_1$ such that $W_z=V$ and $W_1\cdot G_1=1$. In particular $W_1,\ldots, W_{z-1}$ are non-branching in $\wt Q_1$.

Suppose $\psi$ does not touch the chain $W_1+\ldots+W_{z-1}$. Because $G_1$ meets $W_1$, the latter chain consists of non-branching $(-2)$-curves which are contracted by $\pi_0$ subsequently after $C_1$ and $G_1$. Because $\wt Q_1$ contracts to a smooth point we necessarily have $T=[(2)_t,z+1]$ and $V^2=-t-3$ for some $z,t\geq 1$. In fact $z\geq 2$, otherwise either $\psi$ creates a $(-1)$-curve in $D_1-E_1$ or it touches $G_1$, which was already shown to be impossible. Let $U_1$, $U_2$ be the connected components of $\wt Q_1-V$ not containing $C_1$. One of them, say $U_2$, is not touched by $\psi$, so it consists of $(-2)$-curves, hence it is a chain. Write $U_2=[(2)_{u_2}]$ with $u_2\geq 1$. Then the component of $U_1$ meeting $V$ has self-intersection $-u_2-2\leq -3$, so it is touched or contracted by $\psi$. Recall that $A_0$ meets $T$. Because $\psi(V)^2=-2$,  $\psi$ touches $V$ exactly $t+1\geq 2$ times, so $U_1$ is a chain contracted by $\psi$. We get $U_1=[(2)_{u_1},u_2+2]$. Since $z\geq 2$, there are more than two components in $\wt Q_1$ with self-intersection smaller than $-2$. It follows that $A_0$ meets the tips of $\wt Q_1$ contained in $T$ and $U_1$, hence $u_1=0$. Because $\psi$ touches $V$, it contracts $U_2$, so since $A_0$ meets $\Delta_0^-$we have $t=u_1=0$ and $z=2$. But then $\psi$ touches $G_1$; a contradiction.

Therefore $\psi$ touches $W_1+\ldots+W_{z-1}$, hence $z\geq 2$ and $\psi$ is of type $I$. It follows that $T=[(2)_{m_1},3,(2)_{m_2}]$ for some $m_1\geq 1$, $m_2\geq 0$ and hence that $W_1$ is non-branching with $W_1^2=-m_2-2$. If $m_2\neq 0$ then $A_0$ meets $W_1$ and hence $W_2$ is a $(-2)$-curve, so after the contraction of $C_1$ the image of $\wt Q_1$ contains a chain $[3,(2)_{m_2},1,m_2+2,2]$, which is not negative definite; a contradiction. Thus $m_2=0$, which implies that $W_2^2=-m_1-3\leq -4$ and hence $A_0$ meets $W_2$. Then the image of $\wt Q_1$ after the contraction of $C_1+G_1+W_1+T$ consists of a $(-1)$-curve and some number of $(-2)$-curves. Being contractible to a point, it is of type $[1,2,\ldots,2]$ with $W_z$ as the $(-1)$-tip. But then $V=W_z$ cannot be branching in $\wt Q_1$; a contradiction.

\med Let $V_1$, $V_2$ be the components of $\wt Q_1$ met by $A_0$. Because $\eta=0$, we may assume $V_1$ is a tip of $\wt Q_1$ contained in a $(-2)$-twig of $\wt Q_1$. Denote the maximal twig of $\wt Q_1$ containing $V_1$ by $T_1$. Let $T_2$ be the second twig of $\wt Q_1$ other than $C_1$. Suppose $V_2$ is not a tip of $\wt Q_1$. Then $\psi$ is of type $I$. If $V_2\leq T_1$ then $\psi$ does not touch $T_2$, so $T_2$ consists of $(-2)$-curves, hence $T_1=[(2)_{t_1},t_2+3]$ for some $t_1,t_2\geq 0$. But then $V_1$ and $V_2$ are tips of $T_1$, and then we see that $\psi$ creates a $(-1)$-curve in $D_1-E_1$; a contradiction. Thus $V_2$ is a component of $T_2$ and $\psi$ does not touch $T_2-V_2$. Because $R_1$ consists of $(-2)$-curves, $T_2-V_2$ consists of $(-2)$-curves and $T_1=[(2)_{t_1},3,(2)_{t_2}]$ for some $t_1\geq 1$ and $t_2\geq 0$. Because $V_2$ is not a tip of $\wt Q_1$, $T_2$ has at least two components, so \ref{lem:chains_with_d=1} implies that $T_2=[(2)_{t_3},t_1+3,t_2+2]$ for some $t_3\geq 0$. It follows that $t_2=0$, so after contracting $C_1$ the image of $\wt Q_1$ is of type $[(2)_{t_1},3,1,2,t_1+3,(2)_{t_2}]$. Let $U$ be the component of $\wt Q_1$ meeting $G_1$ which has self-intersection $-3$. The contraction of $D_0+A_0-U$ maps $X_0$ onto $\PP^2$ and $E_0$ onto a cuspidal curve with semi-ordinary cusps. By \ref{lem:E_with_mult<=3_is_rectifiable} the latter is Cremona equivalent to a line; a contradiction.

We are left with the case when both $V_1$, $V_2$ are tips of $\wt Q_1$. Then $f=\psi(\wt Q_1-C_1)$ is a $(-2)$-cycle. Let $U$ be a component of $\wt Q_1-G_1-C_1$ meeting $G_1$. Then $\psi(\wt Q_1-U+A_0)$ is a chain of type $[2,\ldots,2,1]$ and $C_1$ is its tip. Thus the contraction of $\wt Q_1-U+A_0$ maps $X_0$ onto $\PP^2$ and maps $E_0$ onto a cuspidal curve with semi-ordinary cusps; a contradiction.
\end{proof}

\bcor\label{cor:|zeta|<=p2-2} If $\bar E\subseteq \PP^2$ is a complex rational cuspidal curve which is not Cremona equivalent to a line and $(X_n,D_n,E_n)$ is as in \ref{ssec:logMMP} then $$|K_n\cdot(K_n+E_n)|\leq p_2(\PP^2,\bar E)-2\leq 2.$$ \ecor

\begin{proof} By \ref{thm:MAIN1} $p_2(\PP^2,\bar E)\in \{3,4\}$, so $\zeta=K_n\cdot(K_n+E_n)\leq p_2(\PP^2,\bar E)-2$ by \ref{prop:jump>=2}. By \ref{cor:b2_bounds_etc}(iv) $\zeta\geq 2-p_2(\PP^2,\bar E)$.
\end{proof}

\section{The length of the minimalization process}\label{sec:n<=1}

Recall that $n$ is the length of the minimalization process $\psi=\psi_n\circ\ldots\circ\psi_1:(X_0,\frac{1}{2}D_0)\to (X_n,\frac{1}{2}D_n)$ or, equivalently, of the process $\psi'=\psi_n'\circ\ldots\circ\psi_1'\:(X,\frac{1}{2}D)\to (X_n',\frac{1}{2}D_n')$ defined in Section \ref{ssec:logMMP}. Because \ref{thm:MAIN1} gives $p_2(\PP^2,\bar E)\geq 3$, we have now better restrictions on $D_n$. In particular, by \ref{lem:(Xi,Di)_properties}(iv) $n+\ind(D_n')\leq 5-p_2(\PP^2,\bar E)\leq 2$. As we now show, this basically leads to Theorem \ref{thm:MAIN3}.

\begin{proof}[Proof of Theorem \ref{thm:MAIN3}] We have $p_2(\PP^2,\bar E)\geq 3$ by \ref{thm:MAIN1}, so \ref{lem:(Xi,Di)_properties}(iv) gives $\ind(D_n')+n\leq 5-p_2(\PP^2,\bar E)\leq 2$. Suppose $n>1$. Then $n=2$, $p_2(\PP^2,\bar E)=3$ and $\ind(D_2')=0$, so $D_2'$ has no twigs. It follows that $s=0$, i.e. $\psi_0$ does not contract any maximal twig of $D$. In particular, $\tau_j^*=\tau_j-1\geq 1$ for each $j\leq c$, so $\tau^*=\tau-c\geq c$. Let $t(B)$ denote the number of maximal twigs of a divisor $B$. Properties of $A_i$ (see \eqref{eq:Ai}) imply that $t(D_{i+1}')-t(D_i')\in\{0,1,2\}$ for $i<2$, so $t(D)\leq t(D_2')+4=4$ and hence $c\leq 2$. Because of the minimality of the resolution $\pi\:(X,D)\to (\PP^2,\bar E)$, we have $t(D)\geq 3$. By \ref{cor:|zeta|<=p2-2} $|\zeta|\leq 1$. By \ref{prop:(Xi,Di)_properties_for_CN} we have \begin{eqnarray} \label{eq:n=2,I}\gamma_2+\tau^*+(2-\eta)+(n_0-\eta_0) &\leq& 2\zeta+6,\\\label{eq:n=2,II} c+(2-n_{exc})+K_2\cdot R_2 &\leq& 3-\zeta. \end{eqnarray} Since $\gamma_2+\tau^*\geq 4+c\geq 5$, we get $\zeta\in\{0,1\}$. Let $A_0'$ and $A_1'$ be the proper transforms of $A_0$, $A_1$ on $X$. Note that $\psi_1$, which contracts $A_0$, does not touch the proper transform of $A_1$ on $X_0$, so we denote it also by $A_1$. We have $K_2^2=\zeta-K_2\cdot E_2=\zeta+2-\gamma_2$. By \ref{lem:(Xi,Di)_properties}(iii) $\#\psi(\wt Q_1+\wt Q_2)=\rho(X_2)+1=11-K_2^2=9-\zeta+\gamma_2$, so $\#D_2\geq 13$.

\bcl $c=1$. \ecl

Suppose $c=2$. Then $D$ has four maximal twigs and, because $D_2'$ has no twigs, each $A_j'$, $j=0,1$ meets exactly two tips of $D$, so $n_0=n=2$. We have $\eta>0$, otherwise \eqref{eq:n=2,I} gives $5\leq \gamma_2+\tau^*\leq 2\zeta+2\leq 4$, which is impossible. Let $L$ be a component of $\Upsilon_2$, say, $\psi_*^{-1}(L)\leq \wt Q_1$. Because $\Delta_2^+=0$, we have $\beta_{D_2}(L)=2$, so since $D_2$ is snc-minimal, $L$ meets only one component of $D_2-L$. Then one of the $A_j$'s, say $A_0$, meets $\wt Q_1$ twice. Then $A_0\cdot \wt Q_2=0$ and $A_1\cdot \wt Q_1=0$, so $\pi_0(A_0)$ and $\pi_0(A_1)$ are disjoint curves on $\PP^2$; a contradiction.

\med Now $E$ is a maximal twig of $D=D_0'$ and $\tau^*=\tau-1\geq 1$. Since $t(D_2')=0$ we have $n_1\geq 1$, so $n_0\leq 1$. Denote the component of $\wt Q_1-C_1$ meeting $E_0$ by $U$. We have $U\cdot E_0=1$.

\bcl $\eta=0$.\ecl

Suppose $\eta\neq 0$. Let $L$ be a component of $\Upsilon_2\cap \cal L$. Because $A_0\cdot A_1=0$, changing the order of contractions if necessary, we may assume that $\psi_1$ creates the $(-1)$-curve $L$. As noted above $\beta_{D_2}(L)=2$ and $L$ meets some component $B$ of $D_2-L$ twice. It follows that $\psi_*^{-1}(L)\neq U$ and $A_0\cdot E_0=0$. We see that $\eta=\eta_0\geq 1$ and, since $\eta_0\leq n_0\leq 1$, that $\Upsilon_2\cap \cal L$ has only one component. We obtain $n_0=n_1=1$ and $\eta=\eta_0=1$. Denote the proper transforms of $L$ and $B$ on $X_0$ by $L'$ and $B'$ respectively. By \eqref{eq:n=2,II} $K_2\cdot R_2\leq 2$, so $B^2\geq -4$. By \ref{lem:(Xi,Di)_properties}(vi) and \ref{lem:properties_of_(Xi,Di)_for_CN}(i) $B'\cdot E_0+1\geq B\cdot E_2\geq 8+2B^2$. If $B=C_1'$ then we get $\tau_1^*\geq 4$, which is impossible by \eqref{eq:n=2,I}. Thus $B\neq C_1'$ and hence $B'\cdot E_0\leq 1$, so $B^2\in \{-3,-4\}$. Suppose $B\cdot E_2>1$. Then $B'=U$ and $U$ meets a $(-2)$-twig contracted by $\psi_2$. Since $L$ meets $B$ twice, $A_0$ meets some other twig meeting $B'$. But because $t(D_2')=0$, $A_0$ meets also a twig in the connected component of $\wt Q_1-B'$ containing $C_1$, so it cannot create a $(-1)$-curve meeting $B$ twice (note $C_1$ and $U$ are not contracted by $\psi$); a contradiction. Thus $B\cdot E_2\leq 1$, so $B^2=-4$ and $K_2\cdot (R_2-B)=0$. By \eqref{eq:n=2,I} and \eqref{eq:n=2,II} $\zeta=0$, $n_{exc}=2$, $\gamma_2=4$ and $\tau^*=1$.

Suppose $t(D)=3$. Then $\wt Q_1$ is a chain. Since $L$ meets $B$ twice, $A_0$ meets only one connected component of $\wt Q_1-C_1$, the one containing $B'$. Because $A_1\cdot E_1=1$ and $n=n_{exc}$, we infer that the second connected component consists of $(-2)$-curves. Hence $\wt Q_1=[(2)_{t_1},1,t_1+2,(t_2)]$ for some $t_1,t_2\geq 1$. Then $B'=[t_1+2]$ and $A_0$ meets $B'$ and the tip of $\wt Q_1$ contained in $[(2)_{t_2}]$. We see that $\#D_2=5$; a contradiction.

Thus $t(D)=4$, so $\wt Q_1$ is a fork. Because $t(D_2')=0$, now $A_0$ meets two tips of $D_0$ and both maximal twigs $T_1$, $T_2$ of $\wt Q_1$ containing them meet $B'$. The curve $B'$ is the unique branching component of $\wt Q_1$, $L'$ is a component of $T_1+T_2$ and $\psi_1$ contracts $A_0+T_1+T_2-L'$. Put $T_3=\wt Q_1-T_1-T_2-B'$. The curve $A_1$ meets $E_0$ and the tip of $\wt Q_1$ contained in $T_3$. Hence $\psi_2$ does not touch $\psi_1(B')$. Since $\psi_2$ is of type $I$ and creates a (non-superfluous) $(-1)$-curve in $D_2-E_2$, we infer that the connected component of $\wt Q_1-C_1$ met by $A_1$ is a chain of $(-2)$-curves containing $U$. Then $T_3=[(2)_{t_1},1,t_1+2,(2)_{t_2}]$ for some $t_1\geq 1$ and $t_2\geq 0$ or $T_3=[(2)_{t_1},1]$ for some $t_1\geq 0$. The divisor $\wt Q_1$ contracts to a smooth point, so after the contraction of $T_3$ the curve $B'$ becomes a $(-1)$-curve. The inequality $(B')^2\leq B^2\leq -4$, shows that the contraction of $T_3$ touches $B'$ more than once. The latter is possible only if $T_3=[(2)_{t_1},1]$. Because $\psi_2$ contracts $A_1+[(2)_{t_1-1}]$, we see that $\#D_2=5$; a contradiction.

\bcl $\zeta=n_0=\tau^*=1$, $\gamma_2=4$ and $\#\cal C_+=0$. \ecl

The inequality \eqref{eq:n=2,I} gives $\gamma_2+\tau^*+n_0\leq 4+2\zeta$, so $\zeta=1$ and $\tau^*\leq 2$. By \eqref{eq:Kn(Kn+Dn)} $K_2\cdot \psi(\wt Q_1)=-\tau^*-\zeta=-\tau^*-1$. Recall that $n_0\leq 1$. Suppose $n_0=0$. Then $t(D)=3$ and both $A_0$, $A_1$ meet $E_0$, hence $\psi$ contracts exactly $A_0$, $A_1$ and two maximal $(-2)$-twigs of $D_0$ meeting them. It follows that $\wt Q_1$ is a chain and $K_0\cdot \wt Q_1=K_2\cdot \psi(\wt Q_1)+2=1-\tau^*\leq 0$. By \ref{lem:Q_with_small_KQ} $\wt Q_1=[(2)_t,1]$ or $\wt Q_1=[(2)_t,3,1,2]$ for some $t\geq 0$. Then $\#D_2\leq 4$; a contradiction. We obtain $\gamma_2+\tau^*\leq 5$, hence $\gamma_2=4$ and $\tau^*=1$. By \ref{lem:properties_of_(Xi,Di)_for_CN}(vi) $\#\cal C_+=0$.

\med The above claim gives $K_2\cdot \psi(\wt Q_1-C_1)=-1$, or equivalently, $K_2\cdot R_2=n_{exc}-1\leq 1$, so $n_{exc}\geq 1$. We may assume $A_1$ (contracted by $\psi_2$) meets $E_0$ and $A_0$ (contracted by $\psi_1$) does not.

\bcl $\psi_1$ does not create a $(-1)$-curve in $D_1-E_1$. \ecl

Suppose $\psi_1$ creates a $(-1)$-curve in $D_1-E_1$. Denote the image of this $(-1)$-curve on $X_2$ by $L$ and its proper transform on $X_0$ by $L'$.  Since $A_0\cdot E_0=0$, we have $L\cdot E_2=L'\cdot E_0\leq 1$. But if $L\cdot E_2=1$ then the contraction of $L$ maps $E_2$ onto a $(-3)$-curve on a smooth rational surface, which is impossible (see \ref{lem:cuspidal_of_gt_and_khalf}(iii)). Thus $L\cdot E_2=0$. Also, $L\cdot C_1'=0$, because otherwise \ref{lem:properties_of_(Xi,Di)_for_CN}(i) gives $\tau_1\geq 4$, which is false.

Suppose there is a component of $D_2-L$ which meets $L$ at least twice. Clearly, $M\cdot E_2\leq \psi_*^{-1}(M)\cdot E_0+1\leq 2$. By \ref{lem:properties_of_(Xi,Di)_for_CN}(i) $M\cdot E_2+2K_2\cdot M\geq 4$, so since $K_2\cdot R_2=n_{exc}-1\leq 1$, we get $M\cdot E_2=2$, $K_2\cdot M=1$ and $n_{exc}=2$. Then $\psi_*^{-1}(M)=U$ and $\psi_2$ touches $\psi_1(U)$. Since $n_{exc}=2$, $\psi_2$ creates a $(-1)$-curve in $D_2-E_2$. Because $\psi_1$ is of type $I$, this $(-1)$-curve is necessarily $M$, hence $K_2\cdot M=-1$; a contradiction.

The divisor $\psi(\wt Q_1)$ contains a unique reduced cycle of rational curves. The cycle contains $L$, so we can write it as $V+L$, where $L$ is not a component of $V$. Because $L$ meets each component of $D_2$ at most once, $L$ meets two components $V_1$, $V_2$ of $V$. By \ref{lem:properties_of_(Xi,Di)_for_CN}(iii) they are components of $R_2$. We have $R_2\cdot E_2=\psi_*^{-1}(R_2)\cdot E_0+1\leq U\cdot E_1+1\leq 2$. If $V_1^2=V_2^2=-2$ then $f=V_1+2L+V_2$ is nef, hence $0\leq f\cdot (2K_2+E_2)=(V_1+V_2)\cdot E_2-4\leq R_2\cdot E_2-4<0$; a contradiction. Thus $V_1^2=-3$ and $V-V_1$ consists of $(-1)$- and $(-2)$-curves. Write $f_R$ for the sum of components of $f$ which are also components of $R_2$. We have $C_1'\cdot (2K_2+E_2)=0$ and $L'\cdot (2K_2+E_2)\leq \psi_*^{-1}(L')\cdot E_2-1\leq 0$ for every $L'$ in $\cal L$. Therefore for $f=V+2L$ we get $f\cdot (2K_2+E_2)\leq f_R\cdot (2K_2+E_2)-4\leq 2R_2\cdot K_2+R_2\cdot E_2-4\leq 0$. But $f$ is nef, so $f\cdot (2K_2+E_2)\geq 0$. In particular, $f_R\cdot K_2=1$, so $n_{exc}=2$. Also, $f_R\cdot E_2=2$, so $f_R$ contains $\psi(U)$ and $\psi_2$ contracts a $(-2)$-twig meeting $U$. Then $U$ is branching in $\wt Q_1$, and hence it is the unique branching component of $\wt Q_1$. Since $n=n_{exc}$, $\psi_2$ creates a $(-1)$-curve in $D_2-E_2$, which is necessarily $\psi(U)$, hence $\psi_1(U)$ is not a $(-1)$-curve. Because $A_0$ joins tips of a twig of $\wt Q_1$ containing $C_1$ with some other twig meeting $U$, the only (non-superfluous) $(-1)$-curve in $D_1$ it can create is $\psi_1(U)$. But $\psi_2$ touches $\psi_1(U)$, so $L=\psi(U)$ has a non-negative self-intersection; a contradiction with \ref{lem:properties_of_(Xi,Di)_for_CN}(v).

\med The above claim gives $n_{exc}=1$ and hence $K_2\cdot R_2=0$. Let $L\leq \cal L$ be the $(-1)$-curve created by $\psi_2$. Clearly, $L$ meets $E_2$. In fact, since $\gamma_2=4$, by the argument above $L\cdot E_2\geq 2$. By \ref{lem:properties_of_(Xi,Di)_for_CN}(vii) $\psi_*^{-1}(L)\cdot E_0\geq 1$, so $\psi_*^{-1}(L)=U$ and $\psi(U)\cdot E_2=2$, so again $\psi_2$ contracts some $(-2)$-twig $[(2)_t]$, $t\geq 1$, meeting $U$. If $U$ is a branching component of $\wt Q_1$ then the maximal twig of $\wt Q_1$ containing $C_1$ is of type $[(2)_{v},1]$ for some $v\geq 1$ and, since $t(D_2')=0$, $A_0$ meets its tip, so $\psi_1$ touches $C_1$, which is in contradiction with $\cal C_+=0$. Put $u=-U^2-2\geq 0$. If $\wt Q_1$ is a chain then it is of type $[(2)_{t},u+2,1,(2)_u,t+3,(v)]$ for some $v\geq 1$, so since $K_2\cdot R_2=0$, the curve with self-intersection $-t-3\leq -4$ meets $A_0$, and hence $\psi_1$ creates a $(-1)$-curve, in contradiction to the last claim.

We infer that $\wt Q_1$ is a fork and its unique branching component $B$ is different than $U$. Since $\psi_2$ touches $B$, it follows that $U$ is contained in a maximal twig of $\wt Q_1$. Also, $\psi_1$ does not touch $U$, so $U^2=-2$. If $C_1$ does not meet $B$ then the maximal twig containing $U$ is necessarily of type $[(2)_{t+1},1,t+3,(2)_{v}]$, $v\geq 0$ and $[t+3,(2)_{v}]$ is not touched by $\psi$, contradicting the equality $K_2\cdot R_2=0$. Therefore, $C_1$ meets $B$ and the maximal twig $T$ of $\wt Q_1$ containing $U$ is of type $[(2)_{t+1},1]$. Since $\psi(B)$ is smooth, one of the components of $\wt Q_1-T$ meeting $B$, call it $U'$, is not contracted by $\psi$. Then the contraction of $\wt Q_1+A_0-U'+A_1-U$ maps $X_0$ onto $\PP^2$ and $E_0$ onto a unicuspidal curve with a semi-ordinary cusp. By \ref{lem:E_with_mult<=3_is_rectifiable} this is a contradiction.
\end{proof}

\section{The Coolidge-Nagata conjecture}\label{sec:CN}

With the knowledge about the minimalization process from Section \ref{sec:logMMP} and Theorems \ref{thm:MAIN1} and \ref{thm:MAIN3} in hand we can now prove Theorem \ref{thm:MAIN2}. We make use of the log BMY inequality \eqref{eq:diamond} and of the non-negativeness of the intersection of $2K_i+E_i$ with the positive part of the Zariski decomposition of $K_i+D_i$ \eqref{eq:*}.

We keep the assumption that the rational cuspidal curve $\bar E\subseteq \PP^2$ is not Cremona equivalent to a line, i.e.\ we assume $\bar E\subseteq \PP^2$ violates the Coolidge-Nagata conjecture. Let $(X_i,D_i)$ and $(X_i',D_i')$, $i=0,\ldots,n$ be the models defined in \ref{ssec:logMMP}. Recall that $(X_0',D_0')=(X,D)$. By \ref{lem:(Xi,Di)_properties}(i) $D_i'$ is snc-minimal. Recall that the surface $\PP^2\setminus\bar E$ is of log general type and contains no $\C^1$. It follows that $X_i'\setminus D_i'$, which is its open subset, is of log general type and contains no $\C^1$. Put $K_i'=K_{X_i'}$ and $\cal P_i=(K_i'+D_i')^+$. By \ref{lem:Bk} $(K_i'+D_i')^-=\Bk D_i'$. Let $t(D_i')$ be the number of maximal twigs of $D_i'$. Let $E_i'$ and $A_i'$ be the proper transforms of $\bar E$ and respectively of $A_i$ on $X_i'$. By \ref{lem:cuspidal_of_gt_and_khalf}(iii) $-(E_i')^2\geq 4$. Recall that $n_1$ is the number of $A_i$'s contracted by $\psi$ meeting $E_i$, or equivalently, the number of $A_i'$'s meeting $E_i'$. We extend this notation by defining $n_1(j)$ to be the number of $A_i'$'s with $i<j$ meeting $E_i'$. In particular, $n_1(i)\leq i$, $n_1(n)=n_1$ and $n_1(0)=0$. Note that $E_i'$ is a maximal twig of $D_i'$ if and only if $c=1$ and $n_1(i)=0$. Put $t_E(D_i')=1$ if $(c,n_1(i))=(1,0)$ and $t_E(D_i')=0$ otherwise.

Taking a smooth fiber of any $\PP^1$-fibration of $X$ we see that $f\cdot (mK+E)=-2m+f\cdot E$, so $f\cdot (mK+E)\to -\8$ for $m\to \8$. Thus there exists a maximal integer $m_E$, for which $h^0(m_EK+E)\neq 0$. We have $h^0(2K+E)\neq 0$, so $m_E\geq 2$.

\bprop\label{prop:CN_twig_bounds} Let $m_E\geq 2$ be as above. For each $i\leq n$ \begin{equation}\cal P_i^2=p_2(\PP^2,\bar E)+i-2+\ind(D_i')\leq 3 \label{eq:diamond}\end{equation} and \begin{equation} t(D_i')+i\leq p_2(\PP^2,\bar E)+\ind(D_i')+\delta(D_i')+\frac{1}{m_E}\max(n_1(i)+c-2,0).\label{eq:*}\end{equation}
\eprop

\begin{proof} The first relation is just \ref{lem:(Xi,Di)_properties}(iv). Fix $i\leq n$ and denote the sum of maximal twigs of $D_i'$ by $T$. The divisor $G=D_i'-T$ (which we call the \emph{core of $D_i'$}) has the same arithmetic genus as $D_i'$, i.e.\ $p_a(G)=p_a(D_i')=i$, hence $(K_i'+G)\cdot G=2(i-1)$. The maximal twigs are supported on $\Supp \Bk D_i'=\Supp (K_i'+D_i')^-$ and hence they intersect $\cal P_i$ trivially. We get $$\cal P_i\cdot D_i'=\cal P_i\cdot G=(K_i'+D_i'-\Bk D_i')\cdot G= (K_i'+G)\cdot G+\beta_{D_i'}(G)-\Bk D_i'\cdot G.$$ By \ref{lem:Bk}(ii) the coefficient in $\Bk D_i'$ of the component of the twig $T_j$ which meets $G$ is $\delta(T_j)$. Then $\cal P_i\cdot D_i'=2i-2+t(D_i')-\delta(D_i')$. We have $$\cal P_i\cdot E_i'=(K_i'+D_i'-\Bk D_i')\cdot E_i'=(K_i'+E_i')\cdot E_i'+\beta_{D_i'}(E_i')-\Bk D_i'\cdot E_i'.$$ We check easily that $\Bk D_i'\cdot E_i'=-t_E(X_i')$, hence $$\cal P_i\cdot E_i'=-2+n_1(i)+c+t_E(X_i')=\max(n_1(i)+c-2,0).$$
Because $m_EK+E\geq 0$, we have $m_EK_i'+E_i'\geq 0$, so since $\cal P_i$ is nef, we obtain $0\leq \cal P_i\cdot (K_i'+\frac{1}{m_E}E_i')=\cal P_i\cdot(K_i'+D_i')-\cal P_i\cdot D_i'+\frac{1}{m_E}\cal P_i\cdot E_i'$. But $$\cal P_i\cdot(K_i'+D_i')=\cal P_i^2=p_2(\PP^2,\bar E)+i-2+\ind (D_i'),$$ which leads to \eqref{eq:*}.
\end{proof}

\bcor\label{cor:CN_twig_bounds-1} For each $i\leq n$ we have $$ t(D_i')+\frac{5}{2}i+p_2(\PP^2,\bar E)\leq 9+\frac{1}{2}\max(n_1(i)-i+c,2-i)\leq 9+\frac{1}{2}\max (c,2-i).$$ Furthermore, if the first inequality is an equality then $\cal P_i^2=3$ and $\delta(D_i')=\ind(D_i')=5-i-p_2(\PP^2,\bar E)$, in particular all maximal twigs of $D_i'$ have length $1$.
\ecor

\begin{proof} We have $\delta(D_i')\leq \ind (D_i')$, so \ref{prop:CN_twig_bounds} gives $\ind (D_i')\leq 5-i-p_2(\PP^2,\bar E)$ and  \begin{equation} t(D_i')\leq 10-3i-p_2(\PP^2,\bar E)+\frac{1}{m_E}\max(n_1(i)+c-2,0). \end{equation} Since $m_E\geq 2$, the first inequality follows. If the equality holds then $\cal P_i^2=3$ and $\delta(D_i')=\ind (D_i')$, which implies that all maximal twigs of $D_i'$ have length $1$. The second inequality follows from $n_1(i)\leq i$.
\end{proof}

\bcor\label{cor:CN_twig_bounds-2} With the above notation and assumptions $2c+t_E(D)\leq t(D)\leq 7$. Furthermore, if $t(D)=7$ then $n=0$ and $p_2(\PP^2,\bar E)=c=3$. \ecor

\begin{proof} Because the log resolution $\pi\:(X,D)\to (\PP^2,\bar E)$ is minimal, each $Q_j=\pi^{-1}(q_j)$ contains at least two maximal twigs of $D$, so $2c+t_E(D)\leq t(D)$. By \ref{thm:MAIN1} $p_2(\PP^2,\bar E)\geq 3$. Then \ref{cor:CN_twig_bounds-1} gives $t(D)+p_2(\PP^2,\bar E)\leq 9+\max(\frac{1}{2}c,1)\leq 9+\max(\frac{1}{4}t(D),1)$, hence $t(D)\leq 12-\frac{4}{3}p_2(\PP^2,\bar E)$. Suppose $t(D)\geq 8$. It follows that $p_2(\PP^2,\bar E)=3$, $t(D)=8$ and $c=4$. Furthermore, all twigs of $D$ have length $1$ and $\cal P_0^2=3$. The former implies that $Q_j=[3,1,2]$ for every $j\leq c$, hence $\ind(D)=4\cdot(\frac{1}{2}+\frac{1}{3})=\frac{10}{3}$. By \eqref{eq:diamond} $\ind(D)=5-p_2(\PP^2,\bar E)=2$; a contradiction. Thus $t(D)\leq 7$. We infer that $c\leq 3$.

Consider the case $t(D)=7$. By \ref{cor:CN_twig_bounds-1} $p_2(\PP^2,\bar E)=3$. Suppose $n>0$. We have $t(D_1')\geq t(D)-2$, so taking $i=1$ we get $5\leq t(D_1')\leq 3+\frac{1}{2}(n_1(1)+c)$, hence $c=3$ and $n_1(1)=1$. Then $A_0$ meets $E_0$, so $t(D_1')=6$, and the inequality fails. Thus $n=0$.

Suppose $c=2$. Then for $i=0$ we have an equality in \ref{cor:CN_twig_bounds-1}, so $\delta(D)=\ind(D)=2$ and all maximal twigs of $D$ have length $1$. It follows that the first characteristic pair for each cusp produces a chain $[3,1,2]$, so the maximal twigs of $D$ are of types $[2]$, $[2]$, $[3]$, $[3]$, $[x_1]$, $[x_2]$, $[x_3]$ for some $x_1, x_2,x_3\geq 2$ and $\sum_{i=1}^3 \frac{1}{x_i}=2-2\cdot (\frac{1}{2}+\frac{1}{3})=\frac{1}{3}$. It follows that $x_i\geq 4$, so if $T$ is the sum of all maximal twigs then $K\cdot T=x_1+x_2+x_3-4\geq 8$. By \ref{cor:b2_bounds_etc}(iii),(iv) $K\cdot T\leq K\cdot R_0'=p_2(\PP^2,\bar E)-\zeta+c=5-\zeta\leq 6$; a contradiction.

Suppose $c=1$. As above, $\delta(D)=\ind(D)=2$ and all maximal twigs of $D$ have length $1$. They are of types $[3]$, $[2]$, $[\gamma]$, $[x_1]$, $[x_2]$, $[x_3]$, $[x_4]$ for some $x_1, x_2,x_3, x_4\geq 2$, $\gamma\geq 4$ and $\sum_{i=1}^4\frac{1}{x_i}=2-(\frac{1}{2}+\frac{1}{3})-\frac{1}{\gamma}<\frac{7}{6}$. We have $K\cdot T=\gamma+\sum_{i=1}^4 x_i-9$. Now \ref{cor:b2_bounds_etc}(iii),(iv) give $K\cdot T\leq K\cdot R_0'=4-\zeta\leq 5$, so $\sum_{i=1}^4 x_i\leq 10$. If $x_i\leq 3$ for each $i$ then $\sum_{i=1}^4\frac{1}{x_i}\geq \frac{4}{3}$, which is impossible. Thus, say, $x_1\geq 4$. Then $x_1=4$ and $x_i=2$ for $i=2,3,4$, so $\sum_{i=1}^4\frac{1}{x_i}=\frac{7}{4}$; a contradiction.
\end{proof}

\bcor\label{cor:p2=/=4} If $p_2(\PP^2,\bar E)\neq 3$ then $n=0$ and $\bar E$ has only one cusp.
\ecor

\begin{proof} By \ref{thm:MAIN1} $p_2(\PP^2,\bar E)=4$. By \ref{prop:CN_twig_bounds} $n\leq 1$ and if $n=1$ then $D_1'$ has no twigs. Since $\psi_1'$ contracts at most two twigs, it follows that $n=0$, otherwise $D_1$ would have at most two maximal twigs, which is impossible. We have now $\ind(D)\leq 1$, so since the contribution from each cusp is by \ref{lem:ind_using_pairs} bigger than $\frac{1}{2}$, we obtain $c=1$.
\end{proof}

Recall that $Q_j$, $j=1,\ldots,c$ is the reduced exceptional divisor of $\pi\:(X,D)\to (\PP^2,\bar E)$ over the cusp $q_j$. Because the unique $(-1)$-curve of $Q_j$ is not a tip of $Q_j$, we have $K\cdot Q_j\geq 0$. Recall that $\zeta=K_n\cdot (K_n+E_n)=K_n'\cdot (K_n'+E_n')$.

\bprop\label{prop:case_n=0} If $n=0$ then: \beq\label{eq:n=0,I} \frac{1}{2}(\gamma_0+\tau^*)- p_2(\PP^2,\bar E) \leq \zeta \leq p_2(\PP^2,\bar E)-2,\eeq \beq\label{eq:n=0,II} \sum_{j=1}^c K\cdot Q_j = p_2(\PP^2,\bar E)-\zeta,\eeq \beq\label{eq:n=0,III} \ind(D) \leq 5-p_2(\PP^2,\bar E)\leq 2.\eeq
\eprop

\begin{proof} The first inequality follows from \ref{cor:b2_bounds_etc}(iv) and \ref{cor:|zeta|<=p2-2}. The middle equality follows from \ref{lem:(Xi,Di)_properties}(ii). The last inequality follows from \ref{prop:CN_twig_bounds} and Theorem \ref{thm:MAIN1}.
\end{proof}

We are now ready to prove Theorem \ref{thm:MAIN2}.

\begin{proof}[Proof of Theorem \ref{thm:MAIN2}] By \ref{prop:Fujita_boundary_excluded} we may assume $\PP^2\setminus \bar E$ is of log general type. By \ref{cor:CN_twig_bounds-2} we may assume that $\bar E$ has exactly three cusps and $D$ has six or seven maximal twigs, i.e.\ $c=3$ and $t(D)=6,7$. By \ref{cor:p2=/=4} $p_2(\PP^2,\bar E)=3$.

\bcl If $n\neq 0$ then $n=1$, $\ind(D_1')\leq 1$,  $t(D)=6$ and $s\leq 1$. Moreover, $A_0'$ meets a tip of $Q_1$ and a tip of $Q_2$. \ecl

By \ref{thm:MAIN3} $n=1$. By \ref{cor:CN_twig_bounds-1} $t(D_1')\leq\frac{1}{2}(9+n_1)$. Since $t(D_1')\geq t(D)-1-n_0=t(D)-2+n_1$, we get $t(D)\leq 6+\frac{1}{2}(1-n_1)$. Thus $t(D)=6$, so each $Q_j$ is a chain. If $n_1=1$ then in \ref{cor:CN_twig_bounds-1} we have an equality, so maximal twigs of $D_1'$ have length $1$ and $\delta(D_1')=\ind(D_1')=1$. But if maximal twigs of $D_1'$ have length $1$ then each $Q_j$ is of type $[3,1,2]$, so $\delta(D_1')=\frac{5}{2}$. It follows that $n_1=0$, and hence $t(D_1')=4$. Because $t(D)-t(D_1')=2$, the $(-1)$-curve $A_0'$ meets two tips of $D$. By \ref{prop:CN_twig_bounds} $\frac{3}{2}\leq \delta(D_1')+\ind(D_1')$ and $\ind(D_1')\leq 1$. If $A_0$ meets only one $Q_j$, then by \ref{lem:ind_using_pairs} $\ind(D_1')>\frac{1}{2}+\frac{1}{2}$, which is impossible. Thus $A_0'$ meets a tip of, say, $Q_1$ and a tip of $Q_2$. If $s\geq 2$ then at least two maximal twigs of $D$ are $(-2)$-chains and are not touched by $\psi'$ ($\psi_0$ contracts them). But then $\ind(D_1')>2\cdot \frac{1}{2}$; a contradiction.

\bcl $n=0$. \ecl

Suppose $n=1$. By \ref{lem:(Xi,Di)_properties}(i) $D_1'$ is snc-minimal. The above claim shows that $t(D)=6$, so each $Q_j$, and hence each $\wt Q_j$, is a chain. We now look at $D_0$. Because the components of $D_0-E_0$ meeting $E_0$ do not belong to maximal twigs of $D_0$, they are not contracted by $\psi=\psi_1$. Observe also, that if $\psi(L)$ is a component of $\cal L$ then, because $A_0$ meets tips if $\wt Q_1$ and $\wt Q_2$, we have $\beta_{D_0}(L)=\beta_{D_1}(\psi(L))=3$ and hence $\psi(L)\cdot E_1=L\cdot E_0=1$. Such a $\psi(L)$ does not meet $\Delta_1^+$, so we get $\eta=0$. Because the previous claim gives also $s\leq 1$, we have $\tau^*\geq 2$. Now \ref{prop:(Xi,Di)_properties_for_CN}(i) gives $2\zeta\geq \gamma_1+\tau^*-4\geq 2$, so $\zeta\geq 1$. By \ref{cor:|zeta|<=p2-2} $\zeta=1$. Then $\tau^*=2$, $s=1$ and $\gamma_1=4$.  Since $\tau_j^*\leq 1$ for each $j$, \ref{lem:properties_of_(Xi,Di)_for_CN}(vi) shows that $C_j$ is not touched by $\psi$. Then \ref{prop:(Xi,Di)_properties_for_CN}(ii) reads as $(1-n_{exc})+K_1\cdot R_1=0$, which gives $n_{exc}=1$ and $K_1\cdot R_1=0$. But we already proved that if $\psi_1(L)$ is a component of $\cal L$ then $\psi_1(L)\cdot E_1=1$, so the contraction of $\psi(L)$ maps $E_1$ onto a smooth $(-3)$-curve. This is in contradiction with \ref{lem:cuspidal_of_gt_and_khalf}(iii).

\med Let $\ind_j$ be the contribution to $\ind(D)$ coming from maximal twigs of $D$ contained in $Q_j$. By \ref{lem:ind_using_pairs} $\ind_j>\frac{1}{2}$ (note that if $Q_j$ is branched then $\ind_j$ is strictly bigger than $\ind$ computed for the chain created by the first characteristic pair of $Q_j$). By \ref{prop:case_n=0} $$\ind(D)=\ind_1+\ind_2+\ind_3\leq 2$$ and $$K\cdot Q_1+K\cdot Q_2+K\cdot Q_3=3-\zeta,$$ with $\zeta\geq -1$. Suppose $\zeta=-1$. Then \ref{prop:case_n=0} gives $\tau^*=0$, hence $\tau_j=2$, $s_j=1$ for each $j$. If $Q_j$ is a chain this implies $K\cdot Q_j=0$, so by \ref{lem:Q_with_small_KQ}(ii) $\ind_j\geq \frac{5}{6}$. But at least two $Q_j$'s are chains, hence $\ind(D)>2\cdot \frac{5}{6}+\frac{1}{2}>2$; a contradiction. Thus $\zeta\geq 0$.

Suppose $t(D)=6$. If $Q_j$ is a chain as in \ref{lem:Q_with_small_KQ} then the bigger is $k$ the bigger is $\ind_j$. Using \ref{lem:Q_with_small_KQ} we check that $\ind_j\geq \frac{5}{6}$ if $K\cdot Q_j=0$, $\ind_j\geq \frac{11}{15}$ if $K\cdot Q_j=1$ and $\ind_j\geq \frac{19}{28}$ if $K\cdot Q_j=2$. In any case $\ind_j>\frac{2}{3}$. Since $\ind(D)\leq 2$, it follows that there is a $Q_j$ with $K\cdot Q_j\geq 3$, say it is $Q_3$. But then $\zeta=0$ and $K\cdot Q_1=K\cdot Q_2=0$, so again $\ind(D)\geq 2\cdot \frac{5}{6}+\frac{1}{2}>2$; a contradiction.

Therefore $t(D)=7$. Exactly one $Q_j$ is branched, say it is $Q_3$. Let $Q_3'$ be the chain created by the contraction of the maximal twig of $Q_3$ containing the $(-1)$-curve and let $\ind_3'$ be the contribution to $\ind_3$ coming from the twigs contained in the proper transform of $Q_3'$. We have $K\cdot Q_3\geq K\cdot Q_3'+1\geq 1$, hence $$K\cdot Q_1+K\cdot Q_2+K\cdot Q_3'\leq 2-\zeta\leq 2,$$ so $\max(K\cdot Q_1,K\cdot Q_2,K\cdot Q_3')\leq 2$. By the computations above we get $\ind_j>\frac{2}{3}$ for each $j$, so $\ind(D)>2$; a contradiction.
\end{proof}

\bibliographystyle{amsplain}
\bibliography{bibl}

\providecommand{\bysame}{\leavevmode\hbox to3em{\hrulefill}\thinspace}
\providecommand{\MR}{\relax\ifhmode\unskip\space\fi MR }
\providecommand{\MRhref}[2]{%
  \href{http://www.ams.org/mathscinet-getitem?mr=#1}{#2}
}
\providecommand{\href}[2]{#2}
\begin{thebibliography}{10}

\bibitem{BHPV}
Wolf~P. Barth, Klaus Hulek, Chris A.~M. Peters, and Antonius Van~de Ven,
  \emph{Compact complex surfaces}, second ed., Results in Mathematics and
  Related Areas. 3rd Series. A Series of Modern Surveys in Mathematics, vol.~4,
  Springer-Verlag, Berlin, 2004.

\bibitem{Coolidge}
Julian~Lowell Coolidge, \emph{A treatise on algebraic plane curves}, Dover
  Publications Inc., New York, 1959.

\bibitem{Fujita-on_the_topology}
Takao Fujita, \emph{On the topology of noncomplete algebraic surfaces}, J. Fac.
  Sci. Univ. Tokyo Sect. IA Math. \textbf{29} (1982), no.~3, 503--566.

\bibitem{Kawamata}
Yujiro Kawamata, \emph{Addition formula of logarithmic {K}odaira dimensions for
  morphisms of relative dimension one}, Proceedings of the {I}nternational
  {S}ymposium on {A}lgebraic {G}eometry ({K}yoto {U}niv., {K}yoto, 1977)
  (Tokyo), Kinokuniya Book Store, 1978, pp.~207--217.

\bibitem{Koras-ab_moh}
Mariusz Koras, \emph{On contractible plane curves}, Affine algebraic geometry,
  Osaka Univ. Press, Osaka, 2007, pp.~275--288.

\bibitem{Langer}
Adrian Langer, \emph{Logarithmic orbifold {E}uler numbers of surfaces with
  applications}, Proc. London Math. Soc. (3) \textbf{86} (2003), no.~2,
  358--396.

\bibitem{MaSa-cusp}
Takashi Matsuoka and Fumio Sakai, \emph{The degree of rational cuspidal
  curves}, Math. Ann. \textbf{285} (1989), no.~2, 233--247.

\bibitem{MiTs-lines_on_qhp}
M.~Miyanishi and S.~Tsunoda, \emph{Absence of the affine lines on the homology
  planes of general type}, J. Math. Kyoto Univ. \textbf{32} (1992), no.~3,
  443--450.

\bibitem{Miyan-OpenSurf}
Masayoshi Miyanishi, \emph{Open algebraic surfaces}, CRM Monograph Series,
  vol.~12, American Mathematical Society, Providence, RI, 2001.

\bibitem{Kumar-Murthy}
N.~Mohan~Kumar and M.~Pavaman Murthy, \emph{Curves with negative
  self-intersection on rational surfaces}, J. Math. Kyoto Univ. \textbf{22}
  (1982/83), no.~4, 767--777.

\bibitem{Nagata}
Masayoshi Nagata, \emph{On rational surfaces. {I}. {I}rreducible curves of
  arithmetic genus {$0$}\ or {$1$}}, Mem. Coll. Sci. Univ. Kyoto Ser. A Math.
  \textbf{32} (1960), 351--370.

\bibitem{Palka_CN-arxiv}
Karol Palka, \emph{The {C}oolidge-{N}agata conjecture holds for curves with
  more than four cusps},
  \href{http://arxiv.org/abs/1202.3491}{arXiv:1202.3491}.

\bibitem{Palka-minimal_models}
\bysame, \emph{Cuspidal curves, minimal models and {Z}aidenberg's finiteness
  conjecture}, \href{http://arxiv.org/abs/1405.5346}{arXiv:1405.5346}.

\bibitem{Palka-exceptional}
\bysame, \emph{Exceptional singular {$\bold Q$}-homology planes}, Ann. Inst.
  Fourier (Grenoble) \textbf{61} (2011), no.~2, 745--774,
  \href{http://arxiv.org/abs/0909.0772}{arXiv:0909.0772}.

\bibitem{Russell2}
Peter Russell, \emph{Hamburger-{N}oether expansions and approximate roots of
  polynomials}, Manuscripta Math. \textbf{31} (1980), no.~1-3, 25--95.

\bibitem{Yoshihara}
Hisao Yoshihara, \emph{A note on the existence of some curves}, Algebraic
  geometry and commutative algebra, {V}ol.\ {II}, Kinokuniya, Tokyo, 1988,
  pp.~801--804.

\bibitem{OrZa_cusp}
M.~G. Za{\u\i}denberg and S.~Yu. Orevkov, \emph{On rigid rational cuspidal
  plane curves}, Uspekhi Mat. Nauk \textbf{51} (1996), no.~1(307), 149--150,
  translation in Russian Math. Surveys \textbf{51} (1996), no. 1, 179–-180.

\end{thebibliography}
\end{document}